\documentclass[11pt,oneside]{amsart}
\usepackage{packages}

\begin{document}


    \title{Partial Semigroupoid Actions on Sets}
    
    \author[Petasny]{Rafael Haag Petasny}
    \address{Instituto de Matematica e Estatística, Universidade Federal do Rio Grande do Sul, Av. Bento Gonçalves, 9500, 91.509-900, Porto Alegre, Brazil}
    \email{rafaelpetasny@gmail.com}
    
    \author[Tamusiunas]{Thaísa Tamusiunas$^*$}
    \email{thaisa.tamusiunas@gmail.com}
    \thanks{$^*$ Corresponding author}
    
    \date{}
    \dedicatory{}
    
    
    \begin{abstract} 
        We introduce partial semigroupoid actions on sets and demonstrate that each such action admits a universal globalization. Our construction extends the universal globalization for partial category actions given by P. Nystedt (Lundström) and the tensor product globalization for strong partial semigroup actions given by G. Kudryavtseva and V. Laan, thereby unifying the theory of partial actions for both categories and semigroups.
    \end{abstract}

    
    \maketitle


\section{Introduction}

    Partial group actions on sets were introduced by R. Exel in \cite{exel1998}, where he characterizes a class of $C^\ast$-algebras as crossed products by partial actions. Examples of partial group action on sets can be easily obtained from a restriction of a global group action to a subset of a given set. In \cite[Teorema 4.23]{abadie1999}, F. Abadie proved that any partial group action on a set can indeed be viewed as a restriction of a unique up to isomorphism global one. This global action was later called \emph{universal globalization}, due to its universal property. In  \cite[Theorem 3.4]{kellendonk2004}, J. Kellendonk and M. Lawson proved that universal globalizations are initial objects in a certain category. Partial actions and the globalization problem have been investigated in many other cases, such as partial actions of monoids (\cite{hollings2007}, \cite{kudryavtseva2015} and \cite{schroder2004}), semigroups \cite{kudryavtseva2023}, restriction semigroups \cite{gould2009}, ordered groupoids (\cite{gilbert2005} and \cite{lautenschlaeger2024}), groupoids \cite{bagio2012partial}, inverse categories \cite{alves2024}, categories \cite{nystedt2018}, and inverse semigroupoids \cite{demeneghi}. The reader can also find more information on partial actions in the surveys \cite{batista2017}, \cite{dokuchaev2011} and \cite{dokuchaev2019recent}, and in Exel's book \cite{exel2017book}.

    On the other hand, J. Cuntz and W. Krieger introduced a class of $C^\ast$-algebras generated by partial isometries under relations defined by a 0-1 matrix, as seen in \cite{cuntz1980}. These algebras were initially studied for edge matrices of directed graphs and were generalized by A. Kumjian and D. Pask to higher-rank graphs in \cite{pask2000}. In the latter, it was shown that a higher-rank graph $C^\ast$-algebra is isomorphic to the $C^\ast$-algebra associated to the graph's path category. Due to this algebraic-structure construction and the fact that a 0-1 matrix cannot always be converted made into the edge matrix of a graph, Exel introduced in \cite{exel2011} a novel algebraic concept: semigroupoids and their associated $C^\ast$-algebras.
    
    As an algebraic structure, a semigroupoid is a set equipped with a partially defined associative operation. It is worth mentioning that there is another (more restrictive) notion of semigroupoid due to B. Tilson \cite{tilson1987}. The latter can be seen as a directed graph equipped with a partially defined associative operation, and therefore it contains the classes of categories and semigroups. However, the class of Exel's semigroupoids is strictly larger than the class of Tilson's semigroupoids (also referred to as semicategories or categorical semigroupoids). Several properties of semigroupoids and graphed semigroupoids can be found in the work by L. Cordeiro in \cite{cordeiro2019}.

    In this paper, we are interested in studying the globalization problem for partial semigroupoid actions on sets. Our method generalizes the work of P. Nystedt (Lundström) in \cite{nystedt2018}, whose construction of a universal globalization for partial category actions traces back to F. Abadie's approach in \cite[Theorem 1.1]{abadie1999} and Kellendonk and Lawson's approach in \cite[Theorem 3.4]{kellendonk2004}. We also compare our results with the tensor product globalization for partial semigroup actions constructed by G. Kudryavtseva and V. Laan in \cite{kudryavtseva2023}, concluding that our universal globalization coincides with the tensor product globalization when the partial semigroup action is non-degenerate.

    Partial category and semigroup actions extend the notion of partial monoid actions studied by M. Megrelishvili and L. Schröder \cite{schroder2004}, which in turn generalize the classical partial group actions of Exel \cite{exel1998} and of Kellendonk and Lawson \cite{kellendonk2004}. An alternative line of generalization, due to P. Demeneghi and F. Tasca \cite{demeneghi}, addresses the globalization problem for partial actions of inverse semigroupoids. In the final section of this work, we discuss in detail the differences between the semigroupoid and the inverse semigroupoid settings.

    The paper is organized as follows. Section 2 recalls the definition of semigroupoids and introduces the notion of partial actions, together with the construction of the standard restriction of a global action to a subset and the definition of universal globalization. Section 3 develops the construction of a universal globalization for any given partial action. In Section 4, we discuss the relationship between our construction and those previously established for semigroups and categories. Finally, Section 5 highlights the differences between partial actions of semigroupoids and those of inverse semigroupoids.


\section{Global and Partial Actions}

    \subsection{Semigroupoids.} For our purposes, by a semigroupoid we mean an Exel's semigroupoid, which is a more general notion than that of Tilson's semigroupoid.

     \begin{defi} \label{def:semigroupoid} \cite[Definition 2.1]{exel2011}
        A \textbf{semigroupoid} is a triple $S = (S,S^{(2)},\star)$ such that $S$ is a set, $S^{(2)}$ is a subset of $S \times S$, and $\star \colon S^{(2)} \to S$ is an operation which is associative in the following sense: if $r,s,t \in S$ are such that either
        \begin{itemize}
            \item[(i)] $(s,t) \in S^{(2)}$ and $(t,r) \in S^{(2)}$, or
            \item[(ii)] $(s,t) \in S^{(2)}$ and $(st,r) \in S^{(2)}$, or
            \item[(iii)] $(t,r) \in S^{(2)}$ and $(s,tr) \in S^{(2)}$,
        \end{itemize}
        then all of $(s,t)$, $(t,r)$, $(st,r)$ and $(s,tr)$ lie in $S^{(2)}$, and $(st)r = s(tr)$.
    \end{defi}

    From now on, let $S$ be semigroupoid. Furthermore, for each $s \in S$, we denote
        $$ S^s = \{ t \in S \colon (s,t) \in S^{(2)} \} \quad\text{and}\quad {^sS} = \{ t \in S \colon (t,s) \in S^{(2)} \}. $$
   In other words, $S^s$ is the set of elements in $S$ that are right composable with $s$, while $^sS$ is the set of elements in $S$ that are left composable with $s$. Due to the associativity of $S$, we have $S^t = S^{st}$ and ${^sS} = {^{st}S}$, whenever $(s,t) \in S^{(2)}$. A semigroupoid is called \textbf{categorical} if, for every $s,t \in S$, the sets $S^s$ and $S^t$ are either disjoint or equal. By \cite[Proposition 2.14]{cordeiro2019}, a semigroupoid is categorical if and only if, for all $s,t \in S$, the sets ${^sS}$ and ${^tS}$ are either disjoint or equal.\\

   The class of semigroupoids generalizes two relevant algebraic structures: the semigroups, which can be seen as the semigroupoids whose composition is totally defined; and categories, which can be seen as semigroupoids with local identities. We present the definition of those classes for future reference.

   \begin{defi} \label{def:semigroup}
       A \textbf{semigroup} is pair $(\mathcal{S},\cdot)$ consisting of a non-empty set $\mathcal{S}$ and a binary associative operation $\cdot \colon \mathcal{S} \times \mathcal{S} \to \mathcal{S}$, in the sense that, for every triple $s,t,r \in \mathcal{S}$, we have $s \cdot (t \cdot r) = (s \cdot t) \cdot r$.
   \end{defi}

    By definition, semigroups are precisely semigroupoids such that $S^{(2)} = S \times S$. Furthermore, every semigroup is a categorical semigroupoid, since $\mathcal{S}^s = \mathcal{S}^t$ for every pair $s,t \in \mathcal{S}$.

   \begin{defi} \label{defi:category}
        A \textbf{category} is a quintuple $\mathcal{C} = (\mathcal{C}_0,\mathcal{C}_1,D,R,\circ)$, where $\mathcal{C}_0$ is called the set of \textit{objects} of $\mathcal{C}$, $\mathcal{C}_1$ is called the set of \textit{morphisms} of $\mathcal{C}$, $D,R \colon \mathcal{C}_1 \to \mathcal{C}_0$ are the \textit{domain} and \textit{range} functions and $\circ$ is a partially defined operation $\mathcal{C}_1 \times \mathcal{C}_1 \to \mathcal{C}_1$ called \textit{composition}, satisfying the following properties:
        \begin{itemize}
            \item[$\bullet$] The composition $f \circ g$ is defined if and only if $D(f) = R(g)$, and in this case $D(f \circ g) = D(g)$ and $R(f \circ g) = R(f)$.
            \item[$\bullet$] The composition $\circ$ is associative, in the sense that: if $f \circ g$ and $g \circ h$ are defined, then $(f \circ g) \circ h = f \circ (g \circ h)$;
        \end{itemize}
        A morphism $f \in \mathcal{C}_1$ is called an \textit{identity} if $g \circ f = g$, whenever $D(g) = R(f)$, and $f \circ h = h$, whenever $D(f) = R(h)$. In this case $f$ must satisfy $D(f) = R(f)$.
        \begin{itemize}
            \item[$\bullet$] For each object $e \in \mathcal{C}_0$ there is an identity $1_e \in \mathcal{C}_1$ such that $D(1_e) = e = R(1_e)$.
        \end{itemize}
    \end{defi}

    Notice that, if $1_e \circ 1_f$ is defined, then $1_e = 1_f$ and, therefore, $e = D(1_e) = D(1_f) = f$. Consequently, $e \mapsto 1_e$ is an injection and, hence, we can assume $\mathcal{C}_0 \subseteq \mathcal{C}_1 =: \mathcal{C}$. Furthermore, for every $f \in \mathcal{C}_1$, the composition $1_{R(f)} \circ f$ is defined. Hence, if $f$ is an identity, then $f = 1_{R(f)} \in \mathcal{C}_0$, which shows that $\mathcal{C}_0$ is precisely the set of identities of $\mathcal{C}_1$.

    A category can be seen as a semigroupoid such that $\mathcal{C}^{(2)} = \{(f,g) \colon D(f) = R(g)\} \subseteq \mathcal{C} \times \mathcal{C}$, and such that every element $f \in \mathcal{C}$ has a left identity $R(f)$ and a right identity $D(f)$. Moreover, every category is a categorical semigroupoid. In fact, for every $f \in \mathcal{C}$ we have:
        $$ \mathcal{C}^f = \{g \in \mathcal{C} \colon (f,g) \in \mathcal{C}^{(2)}\} = \{g \in \mathcal{C} \colon R(g) = D(f)\} = \{g \in \mathcal{C} \colon (D(f),g) \in \mathcal{C}^{(2)} \} = \mathcal{C}^{D(f)}. $$
    Therefore, if there is $h \in \mathcal{C}^f \cap \mathcal{C}^g$, then $D(f) = R(h) = D(g)$. Hence, $\mathcal{C}^f = \mathcal{C}^{R(h)} = \mathcal{C}^g$.\\
   
    Since we are working with semigroupoids, it is important to highlight the distinction between the definition proposed by Exel and that introduced by Tilson. A \textbf{directed graph} is a quadruple $(G_0,G_1,d,r)$, where $G_0$ is the set of \textit{vertex}, $G_1$ is the set of \textit{edges} and $d,r \colon G_1 \to G_0$ are the \textit{domain} and \textit{range} functions.
    
    \begin{defi} \label{def:Tilson_semigroupoid} \cite[Appendix B]{tilson1987}
        A \textbf{graphed semigroupoid} (also called \textbf{Tilson semigroupoid}) is a quintuple $(G_0,G_1,d,r,\star)$, where: \begin{itemize} \item[(i)] $(G_0,G_1,d,r)$ is a directed graph; \item[(ii)] $\star$ is a function that, for each pair of arrows $(g,h)$ with $d(g) = r(h)$, it associates an arrow $\star(g,h) := gh$ with $d(gh) = d(h)$ and $r(gh) = r(g)$; \item[(iii)] $\star$ is associative, meaning that if $f, g, h \in G_1$ are such that $d(f) = r(g)$ and $d(g) = r(h)$, then $(fg)h = f(gh)$ holds.\end{itemize}
    \end{defi}
    
    \begin{obs} \label{rel_grap_cat}
        In \cite[Theorem 2.15]{cordeiro2019}, it was demonstrated that the class of graphed semigroupoids corresponds exactly to the class of categorical semigroupoids. Furthermore, in the same work, all possible underlying graph structures for any given categorical semigroupoid were fully characterized.
    \end{obs}

    \begin{exe} \label{exe:markov_1}
        Markov semigroupoids constitute a class of semigroupoids that contains the path categories of directed graphs and some semigroupoids that cannot be realized as graphed semigroupoids. They are constructed in the following way. Let $\Sigma$ be a set and $A = (A(x,y))_{x,y \in \Sigma}$ be a 0-1 matrix. Assume that:
        \begin{itemize} 
            \item[$\bullet$] $S$ is the set of sequences $(s_1,\dots,s_n)$, $n \geq 1$, $s_1,\dots,s_n \in \Sigma$, such that $A(s_i,s_{i+1}) = 1$;
            \item[$\bullet$] $S^{(2)}$ is the set of pairs $((s_1,\dots,s_n),(t_1,\dots,t_m)) \in S \times S$ such that $A(s_n,t_1) = 1$;
            \item[$\bullet$] $\ast$ is given by $(s_1,\dots,s_n) \star (t_1,\dots,t_m) = (s_1,\dots,s_n,t_1,\dots,t_m) \in S$, when defined.
        \end{itemize}
        Then $S = (S,S^{(2)},\star)$ is a semigroupoid, called \textbf{Markov semigroupoid}.
    \end{exe}

    Notice that, for a Markov semigroupoid, if there exist  $x,y \in \Sigma$ such that $A(x,y) = A(y,x) = 1$, then $S$ is a infinite set. If furthermore $A(x,x) = 1$ and $A(y,y) = 0$, then $x \in S^x \cap S^y$ and $y \in S^x \setminus S^y$; therefore, $S$ does not admit a compatible graph structure, because is not categorical. The next example shows that a Markov semigroupoid can be finite and categorical.
    
    \begin{exe} \label{exe:stst}
        Let $\Sigma = \{s,t\}$ and $A$ be the $2 \times 2$ matrix such that $A(s,t) = 1$ and $A(x,y) = 0$, for all $(x,y) \neq (s,t)$. The Markov semigroupoid associated to $(\Sigma,A)$ is the set $S = \{s,t,st\}$, where the only composition defined is $s \star t = st$. This Markov semigroupoid can be endowed with a graph structure $(S_0,S_1,d,r)$, in the following way: $S_0$ consists of three vertex, say $S_0 = \{e_1,e_2,e_3\}$, $S_1 = S$, and the functions $d,r \colon S \to S_0$ are given by $d(t) = d(st) = e_1$, $d(s) = r(t) = e_2$ and $r(s) = r(st) = e_3$. The semigroupoid $S$ associated with this graph structure is shown in the following diagram.
        \begin{center}
            \begin{tikzpicture}
                \tikzstyle{every path}=[draw, ->];

                \node (1) at (0,0) {$\overset{e_1}{\bullet}$};
                \node (2) at (2,0) {$\overset{e_2}{\bullet}$};
                \node (3) at (4,0) {$\overset{e_3}{\bullet}$};

                \path (1) to node[below]{$t$} (2);
                \path (2) to node[below]{$s$} (3);
                \path[bend left] (1) to node[above]{$st$} (3);
            \end{tikzpicture}
        \end{center}
    \end{exe}

    Moreover, it is important to note that Markov semigroupoids are not the only type of semigroupoids that may fail to admit a compatible graph structure, as illustrated in the following example.

    \begin{exe} \label{exe:not_markov_1}
        Let $S = \{s_1,s_2,t_1,t_2,0\}$ with $S^{(2)} = \{(s_1,t_1),(s_1,t_2),(s_2,t_2)\}$, and define the product as $s_i \star t_j = 0$ for all $(s_i,t_j) \in S^{(2)}$. Then, $S$ forms a semigroupoid that does not admit a compatible graph structure, since $S^{s_1} = \{t_1,t_2\}$ and $S^{s_2} = \{t_2\}$, meaning it is not categorical. 
        
        We now show that $S$ is not a Markov semigroupoid either. Indeed, if it were Markov, then $\Sigma$ would be the set of sequences of length one in $S$, that is, the set of elements that cannot be written as as compositions of other elements. In this case, $\Sigma = \{s_1,s_2,t_1,t_2\}$. However, we would then have $s_1 \star t_1 = (s_1,t_1) \neq (s_1,t_2) = s_1 \star t_2$, which contradicts the definition of the product $s_i \star t_j = 0$ for all $(s_i, t_j) \in S^{(2)}$. Hence, $S$ is not a Markov semigroupoid.
    \end{exe}

    \subsection{Partial actions of semigroupoids.} In this subsection, we introduce the notion of partial actions of semigroupoids on sets, discuss the concept of globalization, and provide illustrative examples.

    We begin by defining what a partial action is. From this point on, $X$ denotes a set. Given a set $Y$, a function $f \colon Y \to X$, and a subset $A \subseteq X$, we denote by $f^{-1}(A)$ the inverse image of $A$ by $f$. That is, $f^{-1}(A) = \{ x \in Y \colon f(x) \in A \}$.

    \begin{defi} \label{def:action}
        A \textbf{partial action of $S$ on $X$} is a pair $\alpha = (\{{_sX}\}_{s \in S}, \{\alpha_s\}_{s \in S})$ consisting of a collection $\{{_sX}\}_{s \in S}$ of subsets of $X$ and a collection $\{\alpha_s\}_{s \in S}$ of functions $\alpha_s \colon {_sX} \to X$ such that:
        \begin{enumerate} \Not{P}
            \item If $(s,t) \in S^{(2)}$, then $\alpha_t^{-1}({_sX} \cap \alpha_t({_tX})) = {_tX} \cap {_{st}X}$; \label{pact-1}
            \item If $(s,t) \in S^{(2)}$, then $(\alpha_s \circ \alpha_t)(x) = \alpha_{st}(x)$, for every $x \in {_tX} \cap {_{st}X}$. \label{pact-2}
        \end{enumerate}
        In this case, we denote $\alpha_s({_sX}) = X_s$. A \textbf{global action} of $S$ on $X$ is a partial action $\alpha$ such that:
        \begin{enumerate} \Not{G}
            \item If $S^s = S^t \neq \emptyset$, then ${_sX} = {_tX}$; \label{gact-1}
            \item If $(s,t) \in S^{(2)}$, then ${_tX} = {_{st}X}$. \label{gact-2}
        \end{enumerate}
        Furthermore, we say a partial action is \textbf{non-degenerate} if $X = \bigcup_{s \in S} ({_sX} \cup X_s)$. Otherwise, we say that it is \textbf{degenerate}.
    \end{defi}

    The next two lemmas are handy to determine if a given partial action is global.

    \begin{lema} \label{lema:global_charact}
        Let $\{_sX\}_{s \in S}$ be a collection of subsets of $X$ and $\{\alpha_s\}_{s \in S}$ be a collection of functions $\alpha_s \colon {_sX} \to \alpha_s({_sX}) = X_s$. Then, for each pair $(s,t) \in S^{(2)}$, the following are equivalent:
        \begin{itemize}
            \item[(i)] Conditions \eqref{pact-1}, \eqref{pact-2} and \eqref{gact-2} hold for $(s,t)$.
            \item[(ii)] ${_tX} = {_{st}X}$, $X_t \subseteq {_sX}$, and $\alpha_s \circ \alpha_t = \alpha_{st}$.
        \end{itemize}
        Consequently, $(\{_sX\}_{s \in S},\{\alpha_s\}_{s \in S})$ is a global action if and only if it satisfies (ii), for all $(s,t) \in S^{(2)}$, and \eqref{gact-1}.

        \begin{proof}
            Assume (i). Then ${_tX} = {_{st}X}$ is precisely \eqref{gact-2}. Applying $\alpha_t$ to \eqref{pact-1} we obtain
                $$ {_sX} \cap X_t = \alpha_t({_tX} \cap {_{st}X}) = \alpha_t({_tX}) = X_t. $$
            Therefore, $X_t \subseteq {_sX}$. Now, \eqref{pact-2} is $(\alpha_s \circ \alpha_t)(x) = \alpha_{st}(x)$, for every $x \in {_{st}X} = {_tX}$, thus $\alpha_s \circ \alpha_t = \alpha_{st}$ as functions. Conversely, assume (ii). By hypothesis, \eqref{gact-2} holds, and the condition $\alpha_s \circ \alpha_t = \alpha_{st}$ implies \eqref{pact-2}. Furthermore, we have $\alpha_t^{-1}({_sX} \cap X_t) = \alpha_t^{-1}(X_t) = {_tX} = {_tX} \cap {_{st}X}$. This proves \eqref{pact-1}.

            In conclusion, if (ii) is satisfied for all $(s,t) \in S^{(2)}$, then the pair $(\{_sX\}_{s \in S},\{\alpha_s\}_{s \in S})$ satisfies \eqref{pact-1}, \eqref{pact-2} and \eqref{gact-2}. If it also satisfies \eqref{gact-1}, then it is a global action by definition.
        \end{proof}
    \end{lema}

    \begin{lema} \label{lema:G1implicaG2}
        Let $\alpha$ be a partial action satisfying \eqref{gact-1}. If $S^s \neq \emptyset$, for every $s \in S$, then $(\alpha,X)$ satisfies \eqref{gact-2}. Therefore $\alpha$ is a global action of $S$.

        \begin{proof}
            Let $(s,t) \in S^{(2)}$. Then $S^t = S^{st}$, and consequently $S^t = S^{st} \neq \emptyset$. Therefore \eqref{gact-1} implies ${_tX} = {_{st}X}$.
        \end{proof}
    \end{lema}
    
    \begin{exe} \label{exe:regular_action}
        Let $S$ be a semigroupoid, $X=S$, ${_sX} = S^s$ and $\alpha_s(x) = sx$. Then $\alpha$ is a global action of $S$ on $X$. In fact, if $(s,t) \in S^{(2)}$, then $S^t = S^{st}$ and $tx \in S^s$ if and only if $x \in S^t$, that is, $X_t = {_sX} \cap X_t$. Hence $\alpha_t^{-1}({_sX} \cap X_t) = {_tX} = {_tX} \cap {_{st}X}$. Since $S$ is associative, we have $\alpha_s \circ \alpha_t = \alpha_{st}$.
    \end{exe}
    
    Example \ref{exe:regular_action} illustrates that a semigroupoid always acts globally on itself via left multiplication. When $S$ is specifically a semigroup or a category, for every $t \in S$, there exists some $s \in S$ such that $(s, t) \in S^{(2)}$. In the case of a category, we can choose $s = R(t)$, the left identity of $t$, and for a semigroup, we can choose any $s \in S$. Consequently, $S = \bigcup_{s \in S} S^s = \bigcup_{s \in S} {_sX}$, and the left multiplication action is non-degenerate. This is not always the case for a semigroupoid. For instance, the left multiplication of the semigroupoid $S = \{s, t, st\}$, as defined in Example \ref{exe:stst}, satisfies $\bigcup_{s \in S} ({_sX} \cup X_s) = \{t, st\} \neq S$.

    \begin{exe} \label{exe:not_markov_2}
        Let $S = \{s_1,s_2,t_1,t_2,0\}$ be the semigroupoid of Example \ref{exe:not_markov_1} and $X = \{1,2,3,4\}$. Consider the following subsets of $X$ and the respective functions $\alpha_s \colon {_sX} \to X$:
        \begin{align*}
            {_{s_1}X} &= \{1,2\}, & \alpha_{s_1}(1) &= 2, \ \alpha_{s_1}(2) = 3, &
            {_{s_2}X} &= \{1\}, & \alpha_{s_2}(1) &= 2, \\
            {_{t_2}X} &= \{1,3\}, & \alpha_{t_2}(1) &= 1, \ \alpha_{t_2}(3) = 3, & {_{t_1}X} &= \{2\}, & \alpha_{t_1}(2) &= 1, \\
            {_0X} &= \{1,2\}, & \alpha_0(1) &= 2, \ \alpha_0(2) = 2.
        \end{align*}
        Then $\alpha = (\{_sX\}_{s \in S},\{\alpha_s\}_{s \in S})$ is a partial action of $S$ on $X$. This action is degenerate since $4 \notin \cup_{s \in S}({_sX} \cup X_s)$. Moreover, it is not global because $0 = s_2t_2$ but ${_0X} \neq {_{t_2}X}$.
    \end{exe}

   
    \subsection{Restriction and globalization.} In this section, we present the construction of partial actions by restricting global ones, and motivate our definition of universal globalization.

    \begin{prop} \label{prop:restriction}
        Let $\beta = (\{{_sY}\}_{s \in S}, \{\beta_s\}_{s \in S})$ be a global action of $S$ on a set $Y$ and $X$ be a subset of $Y$. For each $s \in S$, define
            ${_sX} = \{ x \in X \cap {_sY} \colon \beta_s(x) \in X \} = \beta_s^{-1}(X \cap Y_s) \cap X$
        and $\alpha_s = \beta_s|_{_sX} \colon {_sX} \to \beta_s({_sX})$. Then, $\alpha = (\{_sX\}_{s \in S},\{\alpha_s\}_{s \in S})$ is a partial action of $S$ on $X$.

        \begin{proof}
        To prove \eqref{pact-1}, assume for a moment that, for $(s,t) \in S^{(2)}$, the following identities hold:
        \begin{gather*}
            \beta_t^{-1}({_sX} \cap Y_t) \cap X = \alpha_t^{-1}({_sX} \cap X_t), \quad\text{and} \tag{$\ast$} \\
            \beta_s^{-1}(X \cap Y_{st}) \cap Y_t = \beta_s^{-1}(X \cap Y_s) \cap Y_t. \tag{$\ast\ast$}
        \end{gather*}
        Since $\beta$ is a global action, conditions \eqref{pact-1} and \eqref{gact-2} implies $\beta_t^{-1}({_sY} \cap Y_t) = {_{st}Y} = {_tY}$. From \eqref{pact-2} it follows that $\beta_{st} = \beta_s \circ \beta_t$. Therefore, 
        \begin{align*}
            {_{st}X} \cap {_tX} &= \beta_{st}^{-1}(X \cap {Y_{st}}) \cap \beta_t^{-1}(X \cap Y_t) \cap X \\
            &\overset{\eqref{pact-2}}{=} \beta_{t}^{-1}(\beta_s^{-1}(X \cap {Y_{st}})) \cap \beta_t^{-1}(X \cap Y_t) \cap X \\
            &= \beta_t^{-1}{\Big(} \beta_s^{-1}(X \cap Y_{st}) \cap X \cap Y_t {\Big)} \cap X \\
            &\overset{(\ast\ast)}{=} \beta_t^{-1}{\Big(} \beta_s^{-1}(X \cap Y_{s}) \cap X \cap Y_t {\Big)} \cap X \\
            &= \beta_t^{-1}({_sX} \cap Y_t) \cap X \\
            &\overset{(\ast)}{=} \alpha_t^{-1}({_sX} \cap X_t).
        \end{align*}
        This proves that condition \eqref{pact-1} holds. Since $\alpha_s = \beta_s|_{_sX}$, condition \eqref{pact-2} is easily verified.
        
        It remains to prove ($\ast$) and ($\ast\ast$). For the first one, if $x \in \beta_t^{-1}({_sX} \cap Y_t) \cap X \subseteq X \cap {_tY}$, then $\beta_t(x) \in {_sX} \cap Y_t \subseteq X$. Since $x \in {_tX}$, we have $\beta_t(x) = \alpha_t(x) \in X_t$. Therefore,
            $$ x \in \alpha_t^{-1}(\beta_t(x)) \subseteq \alpha_t^{-1}({_sX} \cap X_t). $$
        This shows the inclusion $\beta_t^{-1}({_sX} \cap Y_t) \cap X) \subseteq \alpha_t^{-1}({_sX} \cap X_t)$. For the reverse inclusion, notice that $X_s \subseteq Y_s$. Since $\alpha_t$ is a restriction of $\beta_t$ to a subset of $X$, we obtain
            $$ \alpha_t^{-1}({_sX} \cap X_t) \subseteq \beta_t^{-1}({_sX} \cap X_t) \cap X \subseteq \beta_t^{-1}({_sX} \cap Y_t) \cap X. $$
        This proves ($\ast$). For the second identity, observe that
        \begin{align*}
            Y_{st} = \beta_{st}({_{st}Y}) \overset{\eqref{pact-2}}{=} \beta_s(\beta_t({_{st}Y})) = \beta_s(\beta_t(\beta_t^{-1}({_sY} \cap Y_t))) = \beta_s({_sY} \cap Y_t) \subseteq Y_s.
        \end{align*}
        Since ${_sY} \cap Y_t \subseteq \beta_s^{-1}(\beta_s({_sY} \cap Y_t))$ and $\beta_s^{-1}(X \cap {Y_s}) \subseteq {_sY}$, we obtain
        \begin{align*}
            \beta_s^{-1}(X \cap Y_{st}) \cap Y_t &= \beta_s^{-1}(X \cap Y_s) \cap \beta_s^{-1}(Y_{st}) \cap Y_t \\
            &= \beta_s^{-1}(X \cap Y_s) \cap \beta_s^{-1}(\beta_s({_sY} \cap Y_t)) \cap Y_t \\
            &= \beta_s^{-1}(X \cap Y_s) \cap \beta_s^{-1}(\beta_s({_sY} \cap Y_t)) \cap {_sY} \cap Y_t \\
            &= \beta_s^{-1}(X \cap Y_s) \cap {_sY} \cap Y_t \\
            &= \beta_s^{-1}(X \cap Y_s) \cap Y_t.
        \end{align*}
        This proves ($\ast\ast$).
        \end{proof}
    \end{prop}
    
    \begin{defi}
        The partial action $\alpha$ obtained from the global action $\beta$ in Proposition \ref{prop:restriction} is called the \textbf{restriction} of $\beta$ to the subset $X$. 
    \end{defi}
    
    We wish to define a \textit{globalization} of a partial action $\alpha$ as any global action $\beta$ such that $\alpha$ can be identified with a restriction of $\beta$. To develop this concept, we will recall the standard definitions for the category of partial actions. From now on, a partial action $\alpha$ of $S$ on $X$ will be denoted by $(\alpha, X)$ to emphasize the set on which the action takes place.

    \begin{defi} \label{def:morphism}
        Let $(\alpha,X)$ and $(\beta,Y)$ be partial actions of $S$ and $\varphi \colon X \to Y$ be a function.
        \begin{itemize}
            \item $\varphi$ is called a \textbf{morphism of partial actions} if, for every $s \in S$, we have
                $$ \varphi({_sX}) \subseteq {_sY} \quad\text{and}\quad (\varphi \circ \alpha_s)(x) = (\beta_s \circ \varphi)(x), \ \forall x \in {_sX}. $$

            \item $\varphi$ is said to be an \textbf{embedding} if it is an injective morphism of partial actions such that
                $$ {_sX} = \{ x \in X \colon \varphi(x) \in {_sY}, (\beta_s \circ \varphi)(x) \in \varphi(X) \}, \ \forall s \in S. $$
            Since $\varphi$ is injective, this is equivalent to $\varphi({_sX}) = \beta_s^{-1}(\varphi(X) \cap Y_s) \cap \varphi(X)$, for every $s \in S$.
        \end{itemize}
    \end{defi} 

    In many contexts, a partial action $(\alpha,X)$ is called a \textbf{$S$-act}, and a morphism of partial actions $(\alpha,X) \to (\beta,Y)$ is called a \textbf{$S$-function}. This is the convention, for instance, in \cite{kudryavtseva2023}, \cite{nystedt2018} and \cite{demeneghi}. Here we chose to adopt the terminology of Abadie from \cite{abadie1999} and his subsequential works.

    It follows directly from the definition that if $(\alpha, X)$ is a partial action, then the identity map $id_X \colon X \to X$ is a morphism of partial actions. Furthermore, given two morphisms of partial actions, $\varphi \colon (\alpha, X) \to (\beta, Y)$ and $\psi \colon (\beta, Y) \to (\theta, Z)$, their composition $\psi \circ \varphi \colon (\alpha, X) \to (\theta, Z)$ is also a morphism of partial actions.

    \begin{defi} \label{def:category}
        The \textbf{category of partial actions} of $S$, which we will denote by $\mathcal{A}_p(S)$, is the category whose objects are the partial actions $(\alpha,X)$ of $S$, and whose morphisms are the morphisms of partial actions. The \textbf{category of global actions} of $S$ is the \textit{full subcategory}\footnote{A full subcategory is a subcategory $H$ of the category $G$ such that, if $X,Y \in Ob(H)$, then $Hom_{H}(X,Y) = Hom_G(X,Y)$.} $\mathcal{A}(S)$ of $\mathcal{A}_p(S)$, whose objects are the global actions $(\beta,Y)$.
    \end{defi}

    Recall that, in a category $\mathcal{C}$, a morphism $f \in \mathcal{C}$ is said to be an \textit{isomorphism} if there exists a morphism  $g \in \mathcal{C}$ such that $gf = D(f)$ and $fg = R(f)$. The next result characterize the isomorphisms in the category of partial actions of $S$, and allow us to properly define \textit{globalization}.
   
    \begin{prop} \label{prop:defi_globalization}
	\begin{itemize}
            \item[(i)] The isomorphisms of $\mathcal{A}_p(S)$ are precisely the bijective embeddings.
            \item[(ii)] A partial action $(\alpha,X)$ is isomorphic to a restriction of a global action $(\beta,Y)$ if and only if there exists an embedding $(\alpha,X) \to (\beta,Y)$.
	\end{itemize}
		
	\begin{proof}
            (i) Suppose that $\varphi \colon (\alpha,X) \to (\beta,Y)$ is an isomorphism of $\mathcal{A}_p(S)$. Then there exists a morphism of partial actions $\psi \colon (\beta,Y) \to (\alpha,X)$ such that $\varphi \circ \psi = id_{(\alpha,X)}$ and $\psi \circ \varphi = id_{(\beta,Y)}$. In particular, we have $\psi = \varphi^{-1}$ as a function, thus $\varphi$ is bijective. Now, given $s \in S$, we calculate
		\begin{align*}
			\varphi^{-1}(\beta_s^{-1}(\varphi(X) \cap Y_s) \cap \varphi(X)) = \varphi^{-1}(\beta_s^{-1}(Y \cap Y_s) \cap Y) = \psi({_sY}) \subseteq {_sX}.
		\end{align*}
            For the reverse inclusion, let $x \in {_sX}$. Then $\varphi(x) \in {_sY}$ and $\beta_s(\varphi(x)) = \varphi(\alpha_s(x)) \in \varphi(X)$. Thus, every isomorphism in $\mathcal{A}_p(S)$ is a bijective embedding.
			
            For the converse, suppose that $\varphi \colon (\alpha,X) \to (\beta,Y)$ is a bijective embedding. Then there exists a function $\psi \colon Y \to X$ such that $\varphi \circ \psi = id_Y$ and $\psi \circ \varphi = id_X$. We aim to demonstrate that $\psi$ is a morphism of partial actions. Since $\varphi$ is an embedding, we have
		\begin{align*}
			{_sX} = \varphi^{-1}(\beta_s^{-1}(\varphi(X) \cap Y_s) \cap \varphi(X)) = \varphi^{-1}(\beta_s^{-1}(Y \cap Y_s) \cap Y) = \psi({_sY}).
		\end{align*}
            Therefore, for each $y \in {_sY}$, there exists $x \in {_sX}$ satisfying $x = \psi(y)$, or equivalently, $y = \varphi(x)$. Since $\beta_s \circ \varphi|_{_sX} = \varphi \circ \alpha_s$, it follows that
                $$ (\psi \circ \beta_s)(y) = (\psi \circ \beta_s \circ \varphi)(x) = (\psi \circ \varphi \circ \alpha_s)(x) = (\alpha_s \circ \psi)(y). $$
            This shows that $\psi = \varphi^{-1}$ belongs to $\mathcal{A}_p(S)$. Hence, $\varphi$ is an isomorphism.\\
			
            (ii) Let $(\beta,Y)$ be a global action and $Z \subseteq Y$. Denote by $(\theta,Z)$ the restriction of $\beta$ to $Z$ and suppose that $\varphi \colon (\alpha,X) \to (\theta,Z)$ is an isomorphism (hence, an embedding). It is known that the inclusion $i \colon Z \to Y$ is an embedding $(\theta,Z) \to (\beta,Y)$. We will demonstrate that $i \circ \varphi \colon X \to Y$ is an embedding. For this, observe that
                $$ \varphi^{-1}(\beta_s^{-1}(\varphi(X) \cap Y_s) \cap \varphi(X)) = \{x \in X \colon \varphi(x) \in {_sY}, \beta_s(\varphi(x)) \in \varphi(X)\}. $$
            Let $x \in X$ be such that $(i \circ \varphi)(x) \in {_sY}$ and $\beta_s((i \circ \varphi)(x)) \in (i \circ \varphi)(X)$. Denote $z = \varphi(x)$. Then, the assumption reads $i(z) \in {_sY}$ and $\beta_s(i(z)) \in (i \circ \varphi)(X) \subseteq i(Z)$. Since $i$ is an embedding, it follows that $\varphi(x) \in {_sZ}$. Therefore, $(i \circ \theta_s)(\varphi(x)) = \beta_s((i \circ \varphi)(x)) \in (i \circ \varphi)(X)$. Since $i$ is injective, it follows that $\theta_s(\varphi(x)) \in \varphi(X)$. As $\varphi$ is an embedding, we conclude that $x \in {_sX}$. Thus,
                $$ (i \circ \varphi)^{-1}(\beta_s^{-1}((i \circ \varphi)(X) \cap Y_s) \cap (i \circ \varphi)(X)) \subseteq {_sX}. $$
            The reverse inclusion follows by the same reasoning as the first part of (i). Therefore the composition $i \circ \varphi \colon (\alpha,X) \to (\beta,Y)$ is an embedding.
			
            For the converse, suppose that $(\beta,Y)$ is a global action and $\varphi \colon (\alpha,X) \to (\beta,Y)$ is an embedding. Let $Z = \varphi(X) \subseteq Y$ and $(\theta,Z)$ be the restriction of $\beta$ to $Z$. We will show that $\varphi \colon (\alpha,X) \to (\beta,Y)$ is an isomorphism. By definition, we have ${_sZ} = \{z \in Z \cap {_sY} \colon \beta_s(z) \in Z\}$. From $\varphi({_sX}) \subseteq {_sY} \cap Z$ and $(\beta_s \circ \varphi)({_sX}) = (\varphi \circ \alpha_s)({_sX}) \subseteq Z$, we obtain $\varphi({_sX}) \subseteq {_sZ}$. As $\theta_s = \beta_s|_{_sZ}$, it follows that
                $$ \theta_s \circ \varphi|_{{_sX}} = \beta_s \circ \varphi|_{{_sX}} = \varphi \circ \alpha_s. $$
            Thus, $\varphi \colon (\alpha,X) \to (\beta,Y)$ is a morphism of partial action. Furthermore, we have
                $$ \theta_s^{-1}(\varphi(X) \cap Z_s) = {_sZ} = \beta_s^{-1}(Z \cap Y_s) \cap Z = \beta_s^{-1}(\varphi(X) \cap Y_s) \cap \varphi(X), $$
            and, as $\varphi \colon (\alpha,X) \to (\beta,Y)$ is an embedding, we conclude that
                $$ \varphi^{-1}(\theta_s^{-1}(\varphi(X) \cap Z_s) \cap \varphi(X)) = \varphi^{-1}(\beta_s^{-1}(\varphi(X) \cap Y_s) \cap \varphi(X)) = {_sX}. $$
            Therefore, $\varphi \colon (\alpha,X) \to (\theta,Z)$ is a bijective embedding and, hence, an isomorphism.
	   \end{proof}
    \end{prop}

    \begin{defi} \label{def:globalization}
        Let $(\alpha,X)$ be a partial action of $S$. A \textbf{globalization} of $(\alpha,X)$ is a triple $(\iota,\beta,Y)$ consisting of a global action $(\beta,Y)$ and an embedding $\iota \colon X \to Y$.
    \end{defi}

    The next example shows that a partial action may have many non-isomorphic globalizations.

    \begin{exe} \label{exe:multiple_globalizations}
        Let $\mathcal{S} = \{0,1\}$ be the semigroup with composition $1 \cdot x = x = x \cdot 1$, for any $x \in \mathcal{S}$, and $0 \cdot 0 = 0$. For each $n \in \mathbb{N}$, define
            $$ X^n = \{x_1, \dots, x_n\} = {_0X^n} = {_1X^n}, \quad \alpha_0^n(x_i) = x_1 \quad\text{and}\quad \alpha_1^n(x_i) = x_i, \ \forall i=1,\dots,n. $$
        Then each pair $(\alpha^n,X^n)$ is a global action of $\mathcal{S}$, and each inclusion $i^n \colon X^1 \to X^n$ is an embedding. Hence, each triple $(i^n,\alpha^n,X^n)$ is a globalization of $(\alpha^1,X^1)$. However, there is no bijection $X^i \to X^j$ when $i \neq j$. Therefore $(\alpha^1,X^1)$ has infinitely many non-isomorphic globalizations.
    \end{exe}

    Our interest is to define a class of globalizations in a way that, for each partial action, there exists at most one globalization (up to isomorphism) in that class. However, there are multiple ways to define such class. For instance, in the context of semigroups, \cite[Theorem 3.5]{kudryavtseva2023} and \cite[Theorem 4.7]{kudryavtseva2023} present non-isomorphic globalizations that are unique in their respective class. Here we consider the class of \textit{universal globalizations}, which appear in \cite[Definition 2.24]{abadie1999} as \textit{enveloping actions}.

    \begin{defi} \label{def:universal_globalization}
        Let $(\alpha,X)$ be a partial action. A \textbf{universal globalization} of $(\alpha,X)$ is a globalization $(\iota,\beta,Y)$ of $(\alpha,X)$ such that, for every morphism $\varphi \colon (\alpha,X) \to (\theta,Z)$, where $(\theta,Z)$ is a global action, there exists a unique morphism $\Phi \colon (\beta,Y) \to (\theta,Z)$ such that $\Phi \circ \iota = \varphi$. That is, the following diagram commutes:
            \begin{center}
                \begin{tikzpicture}[scale=2.5]
                    \tikzstyle{every path}=[draw,->];
        
                    \node (X) at (0,0) {$(\alpha,X)$};
                    \node (Y) at (1,0) {$(\beta,Y)$};
                    \node (Z) at (1,-2/3) {$(\theta,Z)$};
        
                    \path (X) to node[above]{$\iota$} (Y);
                    \path (X) to node[below left]{$\varphi$} (Z);
                    \path[dashed] (Y) to node[right]{$\Phi$} (Z);
                \end{tikzpicture}
            \end{center}
    \end{defi}

    The next proposition shows that, with the above definition, a partial action admits at most one (up to isomorphism) universal globalization.

    \begin{prop}
        Let $(\alpha,X)$ be a partial action of $S$ and suppose that $(\iota,\beta,Y)$ and $(\gamma,\theta,Z)$ are universal globalizations of $(\alpha,X)$. Then there is an isomorphism of partial actions $(\beta,Y) \to (\theta,Z)$.

        \begin{proof}
            Since $(\iota,\beta,Y)$ is a universal globalization of $(\alpha,X)$, the function $\iota \colon X \to Y$ is, in particular, a morphism of partial actions. As $(\beta,Y)$ is a global action and $(\gamma,\theta,Z)$ is a universal globalization, there must be a morphism of partial actions $f \colon (\theta,Z) \to (\beta,Y)$ such that $\iota = f \circ \gamma$.

            Interchanging the roles of $(\iota,\beta,Y)$ and $(\gamma,\theta,Z)$ above, we obtain a morphism of partial action $g \colon (\beta,Y) \to (\theta,Z)$ such that $\gamma = g \circ \iota$. Now, notice that $f \circ g \colon (\beta,Y) \to (\beta,Y)$ satisfies
                $$ (f \circ g) \circ \iota = f \circ \gamma = \iota = id_Y \circ \iota. $$
            Since $(\iota,\beta,Y)$ is a universal globalization, the morphism $F \colon (\beta,Y) \to (\beta,Y)$ such that $F \circ \iota = \iota$ must be unique. Therefore $f \circ g = id_Y$. Analogously, the composition $g \circ f \colon (\theta,Z) \to (\theta,Z)$ satisfies
                $$ (g \circ f) \circ \gamma = g \circ \iota = \gamma = id_Z \circ \gamma. $$
            Since $(\gamma,\theta,Z)$ is a universal globalization, it must be $g \circ f = id_Z$. Therefore $f = g^{-1}$. By definition, both $f$ and $g$ are morphisms of partial actions. Thus, $g \colon (\beta,Y) \to (\theta,Z)$ is an isomorphism.
        \end{proof}
    \end{prop}

    It is easy to see that, if $(\alpha,X)$ is a global action, then $(id_X,\alpha,X)$ is the universal globalization of $(\alpha,X)$. In this case, given a morphism $\varphi \colon (\alpha,X) \to (\beta,Y)$, we can chose $\Phi = \varphi$. Therefore, in Example \ref{exe:multiple_globalizations}, $(i^1,\alpha^1,X^1)$ is the unique universal globalization of $(\alpha^1,X^1)$.\\

    Before constructing a universal globalization, we want to point out that, if universal globalizations exist, then they can be characterized as the initial objects of a suitable category. To prove this, let $\mathcal{C}$ and $\mathcal{D}$ be categories. Recall that a covariant functor $F \colon \mathcal{C} \to \mathcal{D}$ is a pair of functions $F_0 \colon \mathcal{C}_0 \to \mathcal{D}_0$ and $F_1 \colon \mathcal{C}_1 \to \mathcal{D}_1$ such that:
    \begin{itemize}
        \item If $f \in \mathcal{C}_1$, then $D(F_1(f)) = F_0(D(f))$ and $R(F_1(f)) = F_0(D(f))$.
        \item If $f \circ g$ is defined in $\mathcal{C}$, then $F_1(f \circ g) = F_1(f) \circ F_1(g)$.
    \end{itemize}
    In this work, all functors are assumed to be covariant, therefore we will omit the term and simply refer to them as \emph{functors}. Now, fix a functor $F \colon \mathcal{C} \to \mathcal{D}$ and $X \in \mathcal{D}_0$. We define a category $(X \downarrow F)$ in the following way:
    \begin{itemize}
        \item $(X \downarrow F)_0$ is the family of morphisms $f \in \mathcal{D}_1$ such that $D(f) = X$ and $R(f) \in F_0(\mathcal{C}_0)$. That is, the objects of $(X \downarrow F)$ are the morphisms $f \colon X \to F_0(Y)$ for some $Y \in \mathcal{C}_0$.

        \item Given $f \colon X \to F_0(Y)$ and $g \colon X \to F_0(Z)$ in $(X \downarrow F)_0$, define $(X \downarrow F)(f,g)$ as the family of morphisms $\varphi \in \mathcal{D}_1(F_0(Y)), F_0(Y'))$ such that $\varphi \circ f = g$. Now, the family of morphisms $(X \downarrow F)_1$ is the union $\bigcup_{f,g \in (X \downarrow F)_0} (X \downarrow F)(f,g)$.

        \item The composition in $(X \downarrow F)$ is defined as the composition in $\mathcal{D}_1$, according to the diagram:

        \begin{center}
            \begin{tikzpicture}
                \tikzstyle{every path}=[draw,->];

                \node (X) at (0,0) {$X$};
                \node (Y1) at (2.5,1.75) {$F_0(Y_1)$};
                \node (Y2) at (2.5,0) {$F_0(Y_2)$};
                \node (Y3) at (2.5,-1.75) {$F_0(Y_3)$};

                \path (X) to node[above left]{$f$} (Y1);
                \path (X) to node[above right]{$g$} (Y2);
                \path (X) to node[below left]{$h$} (Y3);
                \path (Y1) to node[left]{$\varphi$} (Y2);
                \path (Y2) to node[left]{$\psi$} (Y3);
                \path[dashed,bend left] (Y1) to node[right]{$\psi \circ \varphi$} (Y3);

                \node at (7.5,1) {$\varphi \circ f = g,$};
                \node at (7.5,0) {$\psi \circ g = h,$};
                \node at (7.5,-1) {$(\psi \circ \varphi) \circ f = \psi \circ g = h.$};
            \end{tikzpicture}
        \end{center}
    \end{itemize}
    Then $(X \downarrow F) = ((X \downarrow F)_0,(X \downarrow F)_1,D,R,\circ)$ is a category, where $D$ and $R$ are domain and range functions of the category $\mathcal{D}$. This category is called the \textbf{category of objects $F$-under $X$}, and appears in \cite[\textsection 6]{maclane1998} as a particular case of a \textbf{comma category}. In particular, if we fix $(\alpha,X) \in \mathcal{A}_p(S)$ a partial action and $I \colon \mathcal{A}(S) \to \mathcal{A}_p(S)$ the inclusion functor, given by $I_0(\beta,Y) = (\beta,Y)$ for every global action $(\beta,Y)$ of $S$, and $I_1(f) = f$ for every morphism between global actions, then $((\alpha,X) \downarrow I)$ is a category.

    We recall that an \textbf{initial object} in a category $\mathcal{C}$ is an object $X_0 \in \mathcal{C}_0$ such that, for each object $X \in \mathcal{C}_0$, there exists a unique morphism $f \colon X_0 \to X$. In particular, the unique morphism $X_0 \to X_0$ is the identity $1_{X_0}$. Consequently, any two initial objects of $\mathcal{C}$ are isomorphic.

    \begin{prop} \label{prop:universal_objeto_inicial}
        Let $(\alpha,X)$ be a partial action of $S$. Then:
        \begin{itemize}
            \item[(i)] If $(\iota,\beta,Y)$ is a universal globalization for $(\alpha,X)$, then $\iota$ is an initial object in $((\alpha,X) \downarrow I)$.

            \item[(ii)] Suppose that $(\alpha,X)$ admits a universal globalization. If $\gamma \colon (\alpha,X) \to (\theta,Z)$ is an initial object in $((\alpha,X) \downarrow I)$, then $(\gamma,\theta,Z)$ is a universal globalization for $(\alpha,X)$.
        \end{itemize}
        Therefore, if $(\alpha,X)$ has a universal globalization, then the family of universal globalizations of $(\alpha,X)$ is precisely the family of initial objects in the category $((\alpha,X) \downarrow I)$.

        \begin{proof}
            (i) This is straightforward from the definition. An initial object in $((\alpha,X) \downarrow I)$ is a morphism of partial actions $\iota \colon (\alpha,X) \to I(\beta,Y)$ such that, for each morphism $\varphi \colon (\alpha,X) \to I(\theta,Z)$, there is a unique morphism $\Phi \colon (\beta,Y) \to (\theta,Z)$ such that $\Phi \circ \iota = \gamma$. Hence, if $(\iota,\beta,Y)$ is a universal globalization, then $\iota$ is an initial object in $((\alpha,X) \downarrow I)$.\\
            
            (ii) If $\gamma \colon (\alpha,X) \to I(\theta,Z)$ is an initial object in $((\alpha,X) \downarrow I)$, then for each morphism of partial actions $\varphi \colon (\alpha,X) \to I(\theta',Z')$, there is a unique morphism $\Phi \colon (\theta,Z) \to (\theta',Z')$ such that $\Phi \circ \gamma = \gamma'$. Therefore, for $(\gamma,\theta,Z)$ to be a universal globalization, it remains to show that $(\gamma,\theta,Z)$ is, in fact, a globalization for $(\alpha,X)$.

            To prove this, suppose that $(\iota,\beta,Y)$ is a universal globalization of $(\alpha,X)$. From (i) we know that $\iota$ is an initial object in $((\alpha,X) \downarrow I)$. Since initial objects are unique up to isomorphism, there is an isomorphism $\Psi \colon \iota \to \gamma$. That is, $\Psi \colon (\beta,Y) \to (\theta,Z)$ is a bijective embedding such that $\Psi \circ \iota = \gamma$. Since $(\iota,\beta,Y)$ is a globalization, the morphism $\iota$ is an embedding. Since the composition of embeddings is again an embedding, it follows that $\Psi \circ \iota = \gamma \colon (\alpha,X) \to (\theta,Z)$ is an embedding. Therefore, $(\gamma,\theta,Z)$ is a globalization of $(\alpha,X)$ and, hence, a universal globalization.
        \end{proof}
    \end{prop}

   \begin{obs} In the next section, we will prove that any partial action of $S$ admits universal globalization. Then, we obtain from Proposition \ref{prop:universal_objeto_inicial} that universal globalizations of a partial action $(\alpha,X)$ are, in fact, the initial objects of the category $((\alpha,X) \downarrow I)$. This is a generalization of \cite[Corollary 3.6]{kudryavtseva2023} for semigroups and the first part of \cite[Theorem 3.4]{kellendonk2004} for groups. Notice that, in the latter, it is stated that a universal globalization is an initial object in the category \textbf{Global}, which is the full subcategory of $((\alpha,X) \downarrow I)$ whose objects are the embeddings $(\alpha,X) \to I(\beta,Y)$. \end{obs}


\section{The universal globalization} \label{sec:globalization}

    In this section, we prove that every partial semigroupoid action on a set admits a universal globalization, and characterize when the universal globalization of a partial semigroupoid action is non-degenerate.
    
    \begin{defi} \label{defi:rel-1} Let $s, t \in S$. We write $s \sim t$ if and only if one of the following holds:
        \begin{enumerate}\Not{R}
            \item $s=t$;\label{rel-1.1}
            \item $S^s = S^t \neq \emptyset$; \label{rel-1.2}
            \item There exists $u \in {^tS}$ such that $s = ut$; \label{rel-1.3}
        \end{enumerate}
    \end{defi}
    
    The relation $\sim$ is reflexive and transitive. Reflexivity follows from \eqref{rel-1.1}. To see that $s \sim t$ and $t \sim r$ imply $s \sim r$, it is enough to consider that the relations are given by \eqref{rel-1.2} or \eqref{rel-1.3}. If both relations are given by \eqref{rel-1.2}, then $S^s = S^t = S^r \neq \emptyset$. Hence $s \sim r$ by \eqref{rel-1.2}. If $s \sim t$ is given by \eqref{rel-1.2} and $t \sim r$ is given by \eqref{rel-1.3}, then $t = ur$ for some $u \in S$. Therefore, $S^r = S^{ur} = S^t = S^s \neq \emptyset$, and hence $s \sim r$ by \eqref{rel-1.2}. Analogously, if $s \sim t$ is given by \eqref{rel-1.3} and $t \sim r$ is given by \eqref{rel-1.2}, then $s = ut$ for some $u \in S$ and $S^s = S^{ut} = S^t = S^r \neq \emptyset$. Finally, if both relations are given by \eqref{rel-1.3}, then there exist $u,v \in S$ such that $s = ut$ and $t = vr$, concluding that $s = (uv)r$ and $s \sim r$ by \eqref{rel-1.3}.
    
    The relation $\sim$ defined above may not be symmetric. For example, consider the semigroupoid $S = \{s,t,st\}$ given in Example \ref{exe:stst}. Then, $st \sim t$ by \eqref{rel-1.3}. However, $t \not\sim st$ because $t \neq st$, $S^t = \emptyset$, and ${^{st}S} = \emptyset$.

    Denote by $R$ the smallest equivalence relation on $S$ that contains $\sim$. In this case, since $\sim$ is reflexive, the relation $R$ is given by $sRt$ if and only if there are $r_1,\dots,r_n \in S$ such that $s = r_1$, $r_n = t$ and $r_i \sim r_{i+1}$ or $r_{i+1} \sim r_i$ for each $i=1,\dots,n-1$. In this case, we say that there is a chain of length $n$ connecting $s$ to $t$.

    \begin{lema} \label{lema:right_order}
        With the previous notation, we have:
        \begin{itemize}
            \item[(i)] If $(s,t) \in S^{(2)}$, then $st R t$.
            \item[(ii)] If $sRt$, then $S^s = S^t$.
            \item[(iii)] If $(\alpha,X)$ is a global action and $sRt$, then ${_sX} = {_tX}$.
        \end{itemize}

        \begin{proof}
            (i) Let $(s,t) \in S^{(2)}$. Then $s \in {^tS}$ and $st = s \cdot t$, hence, condition \eqref{rel-1.3} implies $st R t$.\\

            (iii) Let $(\alpha,X)$ be any global action of $S$ and suppose that $sRt$. We proceed by induction on the length of the chain connecting $s$ to $t$. If $s=t$, the result is straightforward. If $s \neq t$ and $n=2$, then $sRt$ means that $s \sim t$ or $t \sim s$. If the relation is given by \eqref{rel-1.2}, then \eqref{gact-1} implies that ${_sX} = {_tX}$. If the relation is given by \eqref{rel-1.3}, then there is $u \in S$ such that $s = ut$ or $t = us$. In any case, \eqref{gact-2} implies that ${_sX} = {_{ut}X} = {_tX}$ or ${_tX} = {_{us}X} = {_sX}$.

            Now, assume that $s$ and $t$ are connected by a chain of length $n>2$ and suppose that, whenever $u,v \in S$ are connected by a chain of length $m<n$, then ${_uX} = {_vX}$. Let $r_1,\dots,r_n \in S$ be the chain connecting $s$ to $t$. Then, by the proof of $n=2$, we have ${_sX} = {_{r_2}X}$. Since $r_2,\dots,r_n \in S$ is a chain of length $n-1$ connecting $r_2$ to $t$, it follows from the induction hypothesis that ${_{r_2}X} = {_tX}$. Therefore, ${_sX} = {_{r_2}X} = {_tX}$.\\
            
            (ii) Follows from (iii) applied to the global action provided in Example \ref{exe:regular_action}.
        \end{proof}
    \end{lema}

    When $S$ is a category, the equivalence class of an element $s$ under $R$ corresponds precisely to the set of all elements $t \in S$ such that $D(t) = D(s)$. When $\mathcal{S}$ is a semigroup, the equivalence class of any element $s$ becomes the entire set $S$ (consider semigroups as graphed semigroupoids with only one vertex). Thus, the relation $R$ attempts to emulate the concept of a \textit{domain}, even when the semigroupoid lacks a compatible graph structure. 
    
    We denote the equivalence class of $s$ under $R$ by $\rho_s$. Given a partial action $(\alpha,X)$ of $S$, let
        $$ {_{\rho_s}X} = \bigcup{\Big(} \{_tX \colon tRs\} \cup \{X_t \colon t \in S^s\} {\Big)}. $$
    By Lemma \ref{lema:right_order}(ii), the definition of ${_{\rho_s}X}$ is independent of the choice of representative. From now on, fix $(\alpha,X)$  as a partial action of $S$ and define $S^\delta = S \cup \{\delta\}$, where $\delta$ is a symbol. Let
        $$ D = \{(s,x) \in S \times X \colon x \in {_{\rho_s}X}\} \cup \{(\delta,x) \colon x \in X\} \subseteq S^\delta \times X. $$
    Whenever we write $(a,x) \in D$ with $a=s,t,r$, it is assumed that $a \in S$, otherwise it may be $a \in S^\delta$.
    
    \begin{defi} \label{defi:rel-2}
        Consider $\approx$ to be the smallest equivalence relation on $D$ that contains the relations:
    \begin{enumerate}\Not{}
        \item $(\delta,y) \sim (s,x)$, if $x \in {_sX}$ and $\alpha_s(x) = y$; \label{rel-2.1}
        \item $(s,x) \sim (t,y)$, if there exists $u \in S^t$ such that $s = tu$, $x \in {_uX}$, and $\alpha_u(x) = y$. \label{rel-2.2}
    \end{enumerate}

    More precisely, $(a,x) \approx (b,y)$ if and only if $(a,x) = (b,y)$, or there exist elements $(c_1,z_1),\dots(c_n,z_n) \in D$ such that $(a,x) = (c_1,z_1)$, $(c_n,z_n) = (b,y)$ and, for each $i=1,\dots,n-1$, either $(c_i,z_i) \sim (c_{i+1},z_{i+1})$ or $(c_{i+1},z_{i+1}) \sim (c_i,z_i)$. In this case, we say that there is a \textbf{chain of length $n$} connecting $(a,x)$ to $(b,y)$.
    \end{defi}

    \begin{lema} \label{lema:equivalence_properties}
        With the above notation, the following properties hold:
        \begin{itemize}
            \item[(i)] If $(t,y) \in D$ and $(s,t) \in S^{(2)}$, then $(st,y) \in D$.
            \item[(ii)] If $(s,x) \approx (t,y)$, then $x \in {_sX}$ if and only if $y \in {_tX}$. In this case $\alpha_s(x) = \alpha_t(y)$.
            \item[(iii)] If $(s,x) \approx (t,y)$ and $(p,s),(p,t) \in S^{(2)}$, then $(ps,x) \approx (pt,y)$.
            \item[(iv)] If $(\delta,x) \approx (t,y)$ and $(p,t) \in S^{(2)}$, then $(p,x) \approx (pt,y)$.
            \item[(v)] If $(\delta,x) \approx (\delta,y)$, then $x=y$.
        \end{itemize}

        \begin{proof}
            (i) By Lemma \ref{lema:right_order}(i), $(s,t) \in S^{(2)}$ implies $\rho_{st} = \rho_t$. If $(t,y) \in D$ then $x \in {_{\rho_t}X} = {_{\rho_{st}}X}$. Therefore, $(st,x) \in D$.\\
            
            (ii) If $(s,x) = (t,y)$ there is nothing to prove. Hence, we can assume there is a chain of length $n \geq 2$ connecting $(s,x)$ to $(t,y)$. We proceed by induction, the base cases being $n=2,3$.
            
            If $n=2$, the relation $(s,x) \sim (t,y)$ is given by condition \eqref{rel-2.2}. Therefore $s = tu$ or $t = su$, for some $u \in S$. In the first case it must be $x \in {_uX}$ and $\alpha_u(x)=y$. Thus, if $x \in {_sX} = {_{tu}X}$, then \eqref{pact-1} implies that $x \in {_{tu}X} \cap {_uX} = \alpha_u^{-1}(_tX \cap X_u)$. Hence, $y = \alpha_u(x) \in {_tX}$. On the other hand, if $y \in {_tX}$, then again from \eqref{pact-1}, we have $\alpha_{u}(x) \in {_tX} \cap X_u$, which implies that $x \in \alpha_{u}^{-1}(_tX \cap X_u) \subseteq {_sX}$. In both cases, $\alpha_t(y) = (\alpha_t \circ \alpha_u)(x) = \alpha_s(x)$. The case $t = su$ and $\alpha_u(y) = x$ is analogous.
         
            If $n=3$, then $(s,x) \sim (a,z) \sim (t,y)$, where either $a \in S$ and both relations are given by \eqref{rel-2.2}, or $a = \delta$ and both relations are given by \eqref{rel-2.1}. In the first case, we can apply the argument for $n=2$ twice. In the second case, the claim follows directly, since $(s,x) \sim (\delta,z)$ and $(\delta,z) \sim (t,y)$ imply that $x \in {_sX}$, $y \in {_tX}$, and $\alpha_s(x) = z = \alpha_t(y)$.

            Now, assume $n>3$ and that (ii) holds whenever the chain connecting $(s,x)$ to $(t,y)$ has length $m<n$. Since $n > 3$, we can write the chain as
                $$ (s,x) \sim (a,z) \sim (b,w) \sim \dots \sim (t,y). $$
            If the relation $(s,x) \sim (a,z)$ is given by \eqref{rel-2.2}, then $a \in S$. By the case $n=2$, we have that $x \in {_sX}$ if and only if $z \in {_aX}$, and in this case $\alpha_s(x) = \alpha_a(z)$. Furthermore, since $a \in S$ and the chain connecting $(a,z)$ to $(t,y)$ has length $n-1$, the induction hypothesis gives that $z \in {_aX}$ if and only if $y \in {_tX}$, and in this case $\alpha_a(z) = \alpha_t(y)$. Therefore, $x \in {_sX}$ if and only if $z \in {_aX}$ if and only if $y \in {_tX}$, and in this case $\alpha_s(x) = \alpha_a(w) = \alpha_t(y)$. 
            
            If, on the other hand, the relation $(s,x) \sim (a,z)$ is given by \eqref{rel-2.1}, then $a = \delta$ and the relation $(a,z) \sim (b,w)$ is also given by \eqref{rel-2.1}. Hence $b \in S$. By the case $n=3$, we obtain that $x \in {_sX}$ if and only if $w \in {_bX}$, and in this case $\alpha_s(x) = \alpha_b(w)$. Since $b \in S$ and the chain connecting $(b,w)$ to $(t,y)$ has length $n-2$, the induction hypothesis yields that $w \in {_bX}$ if and only if $y \in {_tX}$, and in this case $\alpha_b(w) = \alpha_t(y)$. Therefore, $x \in {_sX}$ if and only if $w \in {_bX}$ if and only if $y \in {_tX}$, and in this case $\alpha_s(x) = \alpha_b(w) = \alpha_t(y)$. This completes the proof of (ii).\\

            (iii) If $(s,x) = (t,y)$ there is nothing to prove, hence, we assume there is a chain of length $n \geq 2$ connecting $(s,x)$ to $(t,y)$. We proceed by induction, the base cases being $n=2,3$.
            
            If $n=2$, then $(s,x) \sim (t,y)$ is given by \eqref{rel-2.2}. If $s = tu$, then $ps = (pt)u$, therefore, condition \eqref{rel-2.2} implies $(ps,x) \sim (pt,y)$. The case $t = su$ is analogous.
            
            If $n=3$, then $(s,x) \sim (a,z) \sim (t,y)$ where either $a \in S$ and both relation are given by \eqref{rel-2.2} or $a = \delta$ and both relations are given by \eqref{rel-2.1}. In the first case we have $s = au$ or $a = sv$, both which implies ${^sS} = {^aS}$. Thus, the result follows applying the case $n=2$ twice. For the second case, note that $(p,s) \in S^{(2)}$ and $x \in {_sX}$ implies $\alpha_s(x) \in {_{\rho_p}X}$, hence, $(p,\alpha_s(x)) \in D$ and $(ps,x) \sim (p,\alpha_s(x))$ by \eqref{rel-2.2}. On the other hand, if $(s,x) \sim (\delta,z) \sim (t,y)$, then is must be $\alpha_s(x) = z = \alpha_t(y)$. Therefore,
                $$ (ps,x) \overset{\eqref{rel-2.2}}{\sim} (p,\alpha_s(x)) = (p,\alpha_t(y)) \overset{\eqref{rel-2.2}}{\sim} (pt,y). $$

            Now, let $n>3$ and assume that (iii) is valid whenever the chain connecting $(s,x)$ to $(t,y)$ has length $m<n$. Since $n>3$, we can write the chain as
                $$ (s,x) \sim (a,z) \sim (b,w) \sim \dots \sim (t,y). $$
            If the chain starts with relation \eqref{rel-2.2}, then $a \in S$, and the case $n=2$ shows that $(p,a) \in S^{(2)}$. Hence, by the result for $n=2$, we have $(ps,x) \approx (pa,z)$. Since $a \in S$, $(p,a) \in S^{(2)}$ and the chain connecting $(a,z)$ to $(t,y)$ has length $n-1$, the induction hypothesis gives $(pa,z) \approx (pt,y)$. Therefore, $(ps,x) \approx (pa,z) \approx (pt,y)$. 
            
            If, on the other hand, the chain starts with relation \eqref{rel-2.1}, then $x \in {_sX}$. By (ii), we obtain $y \in {_tX}$ and $\alpha_s(x) = \alpha_t(y)$. Hence $(ps,x) \sim (p,\alpha_s(x)) = (p,\alpha_t(y)) \sim (pt,y)$, as in the case $n=3$.\\

            (iv) Any chain from $(\delta,x)$ to $(t,y)$ can be written as $(\delta,x) \sim (r,z) \approx (t,y)$, for some $r \in S$. As $(\delta,x) \sim (r,z)$ implies $z \in {_rX}$ and $\alpha_r(z) = x$, from (ii), we obtain $y \in {_tX}$ and $x = \alpha_r(z) = \alpha_t(y)$. Therefore $(p,x) = (p,\alpha_t(y)) \sim (pt,y)$ by \eqref{rel-2.2}.\\

            (v) If $x = y$ there is nothing to prove. Suppose that $x \neq y$ and $(\delta,x) \approx (\delta,y)$, then the chain connecting $(\delta,x)$ to $(\delta,y)$ must contain the relation \eqref{rel-2.1}. Therefore, it is enough to consider a chain of the form
                $$ (\delta,x) \overset{\eqref{rel-2.1}}{\sim} (t,z) \approx (r,w) \overset{\eqref{rel-2.1}}{\sim} (\delta,y). $$
            In this case, from condition \eqref{rel-2.1} and (ii), we obtain $x = \alpha_t(z) = \alpha_r(w) = y$, a contradiction. Hence, $(\delta,x) \approx (\delta,y)$ only occurs when $x=y$.
        \end{proof}
    \end{lema}

    \begin{defi} \label{def:globalization_action}
        We will denote by $E$ the quotient set $D/\approx$, and the equivalence class of $(a,x) \in D$ by $[a,x] \in E$. For each $s \in S$, let
            $$ {_sD} = \{ (a,x) \in D \colon (a,x) \approx (t,y), \text{ for some } t \in S^s\} \cup \{(\delta,x) \in D \colon x \in {_{\rho_s}X}\}, $$
        and denote by ${_sE}$ the quotient set ${_sD}/\approx$. Furthermore, for each $[a,x] \in {_sE}$, define
        \begin{align*}
            \beta_s([a,x]) = \begin{cases}
                [st,y], &\text{if } (a,x) \approx (t,y) \in {_sD}, \\
                [s,x], &\text{if } a = \delta.
            \end{cases}
        \end{align*}
    \end{defi}

    We now show that $\beta_s \colon {_sE} \to E$ is a well-defined function. Let $[a,x] \in {_sE}$. First, we prove that $\beta_s([a,x]) \in E$. Indeed, if $(a,x) \approx (t,y)$ and $st$ is defined, then $(st,y) \in D$ by Lemma \ref{lema:equivalence_properties}(i). On the other hand, if $(\delta,x) \in {_sD}$, then $x \in {_{\rho_s}X}$. Hence, $(s,x) \in D$ by definition. In both cases, we conclude that $\beta_s([a,x]) \in {D/\approx} = E$. Next, we show that $\beta_s([a,x])$ does not depend on the choice of the representatives. Suppose that $(a,x) \approx (t,y)$ and $(a,x) \approx (r,z)$, with both $st$ and $sr$ defined. Then $[st,y] = [sr,z]$ by Lemma \ref{lema:equivalence_properties}(iii), which implies that $\beta_s([a,x])$ is independent of the representative $(t,y)$. Finally, assume that $a = \delta$ and $(\delta,x) \approx (t,y)$ for some $(t,y) \in D$ with $st$ defined. By Lemma~\ref{lema:equivalence_properties}(iv), we have $[s,x] = [st,y]$, showing that both definitions of $\beta_s([a,x])$ coincide whenever ambiguity could arise. Hence, $\beta_s \colon {_sE} \to E$ is well defined.\\

    The next example illustrates the difference between the sets $S^\delta \times X$, $D$, and $E$. In particular, we show that $E$ may not be equal to the union of the sets ${_sE}$ and $\beta_s({_sE})$. That is, viewed as a partial action, the pair $(\beta,E)$ can be degenerate.

    \begin{exe} \label{exe:not_markov_4}
        Let $S = \{s_1,s_2,t_1,t_2,0\}$ be the semigroupoid of Example \ref{exe:not_markov_1}. Then:
        \begin{gather*}
            \rho_{s_1} = \{s_1\}, \quad \rho_{s_2} = \{s_2\} \quad\text{and}\quad \rho_{t_1} = \{t_1,t_2,0\} = \rho_{t_2} = \rho_0.
        \end{gather*}
        For $X = \{1,2,3,4\}$ and the partial action $(\alpha,X)$ defined in Example \ref{exe:not_markov_2}, we have:
        \begin{gather*}
            {_{\rho_{s_1}}X} = \{1,2,3\} = {_{\rho_{0}}X} \quad\text{and}\quad {_{\rho_{s_2}}X} = \{1,3\}.
        \end{gather*}
        Therefore $S^\delta \times X$ has 24 elements, while $D$ has 18 elements. Namely,
        \begin{gather*}
            (\delta,n), \ n=1,2,3,4, \\
            (s_1,n), (t_1,n), (t_2,n), (0,n), \ n=1,2,3, \\
            (s_2,n), \ n=1,3.
        \end{gather*}
        The equivalence classes of $D$ under $\approx$ are:
        \begin{gather*}
            [\delta,1] = [t_1,2] = [t_2,1], \\
            [\delta,2] = [s_1,1] = [s_2,1] = [0,1] = [0,2], \\
            [\delta,3] = [s_1,2] = [t_2,3], \\
            [s_1,3] = [0,3] = [s_2,3], \\
            [t_1,1], \quad 
            [t_1,3], \quad
            [t_2,2], \quad [\delta,4].
        \end{gather*}
        Therefore $E$ has 8 elements. The sets ${_sE}$ are given by:
        \begin{align*}
            {_{s_1}E} &= \{ [\delta,1], [\delta,2], [\delta,3], [t_1,1], [t_1,3], [t_2,2], [t_2,3] \}, \\
            {_{s_2}E} &= \{ [\delta,1], [\delta,3], [t_2,2], [t_2,3] \}, \\
            {_{t_1}E} &= \{ [\delta,1], [\delta,2], [\delta,3] \} = {_{t_2}E} = {_0E}.
        \end{align*}
        Notice that $[0,3] = [s_1t_1,3] = \beta_{s_1}([t_1,3]) \in E_{s_1}$. However, the class $[\delta,4]$ does not belong to any ${_sE}$ or $E_s$. Therefore $\cup_{s \in S} ({_sE} \cup E_s) = E \setminus \{[\delta,4]\}$.
    \end{exe}

    \begin{teo} \label{teo:globalization}
        Define $\iota \colon X \to E$ by $\iota(x) = [\delta,x]$. Then:
        \begin{itemize}
            \item[(i)] The pair $(\beta,E)$ is a global action of $S$ on $E$.
            \item[(ii)] The triple $(\iota,\beta,E)$ is a universal globalization for $(\alpha,X)$.
        \end{itemize}

        \begin{proof}
            (i) We first verify \eqref{gact-1}. If $s,t \in S$ are such that $S^s = S^t \neq \emptyset$ and $[a,x] \in {_sE}$, then $[a,x] = [r,y]$ for some $r \in S^s$, or $[a,x] = [\delta,y]$ for some $y \in {_{\rho_s}X}$. In the first case we have $r \in S^s = S^t$ and, hence, $[a,x] \in {_tE}$. In the second, $S^s = S^t \neq \emptyset$ implies $\rho_s = \rho_t$, hence $y \in {_{\rho_t}X}$ and $[a,x] = [\delta,y] \in {_tE}$. This shows ${_sE} \subseteq {_tE}$. Since $S^s = S^t \neq \emptyset$ is a symmetric condition, the same argument proves ${_sE} = {_tE}$. For the remaining of the proof for (i) we uses Lemma \ref{lema:global_charact}(ii). That is, for $(s,t) \in S^{(2)}$ we must show ${_tE} = {_{st}E}$, $\beta_t({_tE}) \subseteq {_sE}$ and $\beta_s \circ \beta_t = \beta_{st}$.
            
            If $[a,x] \in {_tE}$, then $\beta_t([a,x]) = [tr,y]$ for some $r \in S^t$ or $\beta_t([a,x]) = [t,y]$. Since $t,tr \in S^s$, in both cases we have $\beta_t([a,x]) \in {_sE}$. Thus $E_t \subseteq {_sE}$. Since $(s,t) \in S^{(2)}$ implies $S^t = S^{st}$ and, by Lemma \ref{lema:right_order}(i), $\rho_t = \rho_{st}$, we have
            \begin{align*}
                {_tE} = \{ [(r,x)] \colon r \in S^{st}\} \cup \{[\delta,x] \colon x \in {_{\rho_{st}}X}\} = {_{st}E}.
            \end{align*}
            Now, given $[a,x] \in {_tE} = {_{st}E}$, we calculate
            \begin{align*}
                (\beta_s \circ \beta_t)([a,x]) = \begin{cases}
                    \beta_s([ta,x]),& \text{if } a \in S^t, \\
                    \beta_s([t,x]),& \text{if } a = \delta.
                \end{cases} = \begin{cases}
                    [sta,x],& \text{if } a \in S^{st}, \\
                    [st,x],& \text{if } a = \delta.
                \end{cases} = \beta_{st}([a,x]).
            \end{align*}
            Therefore $(\beta,E)$ is a global action.\\
            
            (ii) We first verify that $(\iota,\beta,E)$ is a globalization for $(\alpha,X)$. Let $x \in {_sX}$ and notice that ${_sX} \subseteq {_{\rho_s}X}$. Then $(\delta,x) \in {_sD}$ and, therefore, $\iota(x) = [\delta,x] \in {_sE}$. Furthermore, for such $x$ we have
                $$ (\beta_s \circ \iota)(x) = \beta_s([\delta,x]) = [s,x] \overset{\eqref{rel-2.1}}{=} [\delta,\alpha_s(x)] = (\iota \circ \alpha_s)(x). $$
            Thus, $\iota$ is a morphism of partial actions. It follows from Lemma \ref{lema:equivalence_properties}(v) that $\iota$ is injective. Therefore, to $\iota$ to be an embedding, it remains to prove that
                ${_sX} = \{x \in X \colon \iota(x) \in {_sE}, \beta_s(\iota(x)) \in \iota(X)\}$.
            Suppose that $x \in X$ is such that $\iota(x) \in {_sE}$ and $\beta_s(\iota(x)) \in \iota(X)$. From $\iota(x) \in {_sE}$ we have two possible cases:
            \begin{itemize}
                \item[$\bullet$] $[\delta,x] = [\delta,y]$ for some $y \in {_{\rho_s}X}$. In this case, Lemma \ref{lema:equivalence_properties}(v) implies $x=y$. Therefore $\beta_s(\iota(x)) = [s,x]$. Now, if $[s,x] \in \iota(X)$, then there is $z \in X$ such that $[s,x] = [\delta,z]$. As any chain from $(\delta,z)$ to $(s,x)$ must use relation \eqref{rel-2.1}, it follows from Lemma \ref{lema:equivalence_properties}(ii) that $x \in {_sX}$.

                \item[$\bullet$] $[\delta,x] = [t,y]$, for some $t \in S^s$. Similar to the previous case, we obtain from Lemma \ref{lema:equivalence_properties}(ii) that $y \in {_tX}$ and $\alpha_t(y) = x$. Now, if $\beta_s(\iota(x)) = [st,y] \in \iota(X)$, then there is $z \in X$ such that $[st,y] = [\delta,z]$. Again, it must be $y \in {_{st}X}$. Therefore
                    $$ x = \alpha_t(y) \in \alpha_t({_tX} \cap {_{st}X}) \overset{\eqref{pact-1}}{=} \alpha_t(\alpha_t^{-1}({_sX} \cap X_t)) = {_sX} \cap X_t. $$
            \end{itemize}
            This shows that $\{x \in X \colon \iota(x) \in {_sE}, \beta_s(\iota(x)) \in \iota(X)\} \subseteq {_sX}$. The reverse inclusion follows from $\iota$ being a morphism of partial actions. This concludes that $(\iota,\beta,E)$ is a globalization for $(\alpha,X)$.
            
            Now we verify that $(\iota,\beta,E)$ is a universal globalization. To prove this, let $(\theta,Y)$ be a global action and $\varphi \colon (\alpha,X) \to (\theta,Y)$ be a morphism. Define $\Phi \colon D \to Y$ as
            \begin{align*}
                \Phi(\delta,x) = \varphi(x) \quad\text{and}\quad \Phi(s,x) = (\theta_s \circ \varphi)(x).
            \end{align*}
            For $\Phi$ to be well defined, we need to verify that $\varphi(x) \in {_sY}$ whenever $(s,x) \in D$. In fact, if $(s,x) \in D$, then $x \in {_{\rho_s}X}$ and, hence, we have two possible cases:
            \begin{itemize}
                \item[$\bullet$] $x \in {_tX}$ for some $t \in S$ such that $sRt$. In this case, since $\varphi$ is a morphism of partial actions, we have $\varphi(x) \in {_tY}$. As $(\theta,Y)$ is a global action, it follows from Lemma \ref{lema:right_order}(iii) that $\varphi(x) \in {_tY} = {_sY}$.
                
                \item[$\bullet$] $x \in {X_t}$ for some $t \in S$ such that $(s,t) \in S^{(2)}$. In this case we have $x = \alpha_t(y)$, for some $y \in {_tX}$. Since $\varphi$ is a morphism of partial actions we have $\varphi(x) = (\varphi \circ \alpha_t)(y) = (\theta_t \circ \varphi)(y) \in Y_t$. As $(\theta,Y)$ is a global action, it follows from Lemma \ref{lema:global_charact}(ii) that $\varphi(x) \in Y_t \subseteq {_sY}$.
            \end{itemize}
            Hence, the function $\Phi$ is well defined. Furthermore, $\Phi$ is compatible with the relation $\sim$ on $D$. For if $(s,x) \sim (\delta,y)$, then $x \in {_sX}$ and $y = \alpha_s(x)$. Therefore,
                $$ \Phi(s,x) = (\theta_s \circ \varphi)(x) = (\varphi \circ \alpha_s)(x) = \Phi(\delta,\alpha_s(x)) = \Phi(\delta,y). $$
            While if $(s,x) \sim (t,y)$, then there is $u \in S^t$ such that $s = tu$, $x \in {_uX}$ and $\alpha_u(x) = y$. Therefore,
                $$ \Phi(s,x) = (\theta_s \circ \varphi)(x) = (\theta_t \circ \theta_u \circ \varphi)(x) = (\theta_t \circ \varphi)(\alpha_u(x)) = (\theta_t \circ \varphi)(y) = \Phi(t,y). $$
            The case $t=su$ is analogous. As $\Phi$ is compatible with $\sim$ and, consequently, with $\approx$, it induces a function $\Phi \colon E \to Y$, given by $\Phi([a,x]) = \Phi(a,x)$. The function $\Phi$ is a morphism of partial actions and satisfy $\Phi \circ \iota = \varphi$. For if $[t,x] \in {_sE}$, where $t \in S^s$, then
                $$ (\Phi \circ \beta_s)([t,x]) = \Phi([st,x]) = (\theta_{st} \circ \varphi)(x) = (\theta_s \circ \theta_t \circ \varphi)(x) = (\theta_s \circ \Phi)([t,x]), $$
            while if $[\delta,x] \in {_sE}$, then
                $$ (\Phi \circ \beta_s)([\delta,x]) = (\theta_s \circ \varphi)(x) = (\theta_s \circ \Phi)([\delta,x]). $$
            Thus, $\Phi$ is a morphism of partial action. For every $x \in X$, we have $(\Phi \circ \iota)(x) = \Phi([\delta,x]) = \phi(x)$. For the uniqueness of $\Phi$, suppose that $\Psi \colon E \to Y$ is a morphism such that $\Psi \circ \iota = \varphi = \Phi \circ \iota$. For each $[\delta,x] \in E$, we have $\Psi([\delta,x]) = (\Psi \circ \iota)(x) = (\Phi \circ \iota)(x) = \Phi([\delta,x])$. Now, for each $[s,x] \in E$, we have $x \in {_{\rho_s}X}$. Therefore $[s,x] = \beta_s([\delta,x])$ and, hence,
                $$ \Psi([s,x]) = (\Psi \circ \beta_s)([\delta,x]) = (\theta_s \circ \Psi)([\delta,x]) = (\theta_s \circ \Phi)([\delta,x]) = \Phi([s,x]). $$
            That is, $\Psi \equiv \Phi$. Therefore $(\iota,\beta,E)$ is a universal globalization for $(\alpha,X)$.
        \end{proof}
    \end{teo}

    To end this section, we will prove that the universal globalization of a partial semigroupoid action is non-degenerate if and only if the partial action is non-degenerate. Let $(\alpha,X)$ be a partial semigroupoid action. We denote
        $$ X_1 = \cup_{s \in S} ({_sX} \cup X_s) \quad\text{and}\quad X_0 = X \setminus X_1. $$
    That is, $X_1$ and $X_0$ are, respectively, the non-degenerate and the degenerate parts of $X$. With this notation, $(\alpha,X)$ is non-degenerate if and only if $X_1 = X$ if and only if $X_0 = \emptyset$.
    
    \begin{prop} \label{prop:degenerate-1}
        Let $(\alpha,X)$ be a partial semigroupoid action of $S$, and let $(\iota,\beta,E)$ be its universal globalization, as stated in Theorem \ref{teo:globalization}. Then $\iota(X_0) = E_0$. Consequently, $(\alpha,X)$ is non-degenerate if and only if $(\beta,E)$ is non-degenerate.

        \begin{proof}
            First, we prove the inclusion $\iota(X_0) \subseteq E_0$. Let $x \in X_0$ and $\iota(x) = [\delta,x] \in E$. Suppose that $\iota(x) \in {_sE}$, for some $s \in S$, and recall that ${_sE} = {_sD}/\approx$, where
                $$ {_sD} = \{ (a,x) \in D \colon (a,x) \approx (t,y), t \in S^s \} \cup \{ (\delta,x) \colon x \in {_{\rho_s}X}\}. $$
            Since $x \in X_0$, it follows that $x \notin {_{\rho_s}X}$. Therefore,  $(\delta,x) \approx (t,y)$ for some $t \in S^s$. This implies that there exists a chain of the form
                $$ (\delta,x) \overset{\eqref{rel-2.1}}{\sim} (r,z) \approx (t,y), $$
            for some $r \in S$. However, $(\delta,x) \sim (r,z)$ means that $x = \alpha_r(z) \in X_r \subseteq X_1$, a contradiction. Now, suppose that $\iota(x) \in E_s$, for some $s \in S$. Then, it must be
                $$ [\delta,x] = \beta_s([\delta,y]) = [s,y] \quad\text{or}\quad [\delta,x] = \beta_s([t,y]) = [st,y]. $$
            Both cases lead to $(\delta,x) \approx (r,y)$, for some $r \in S$, which, again, result in a contradiction. Therefore,
                $$ \iota(x) \notin \cup_{s \in S} ({_sE} \cup E_s) = E_1 \implies \iota(x) \in E \setminus E_1 = E_0. $$
            This concludes that $\iota(X_0) \subseteq E_0$.
            
            To prove the inclusion $E_0 \subseteq \iota(X_0)$, we will first consider the case $[a,y] \in E_0 \cap \iota(X)$. That is, $[a,y] = [\delta,x]$, for some $x \in X$. Suppose that $x \in X_1$. Then there is $s \in S$ such that $x \in {_sX} \cup X_s$. Since $\iota \colon (\alpha,X) \to (\beta,E)$ is a morphism of partial actions, if $x \in {_sX}$, then
                $$ [a,y] = \iota(x) \in \iota({_sX}) \subseteq {_sE}. $$
            On the other hand, if $x \in X_s$, then
                $$ [a,y] = \iota(x) = (\iota \circ \alpha_s)(z) = (\beta_s \circ \iota)(z) \in E_s. $$
            In both cases, we have a contradiction with $[a,y] \in E_0$. Therefore, $E_0 \cap \iota(X) \subseteq \iota(X_0)$. Now, chose an element $[a,y] \in E \setminus \iota(X)$. Then, we can fix a representative $(s,x) \approx (a,y)$ such that $s \in S$. Since $(s,x) \in D$ and $s \in S$, it follows from the definition of $D$ that
                $$ x \in {_{\rho_s}X} = \bigcup {\Big(} \{ {_tX} \colon tRs \} \cup \{X_t \colon t \in S^s\} {\Big)}. $$
            If $x \in {_tX} \subseteq {_{\rho_s}X}$ then, by the definition of ${_sD}$, it must be
                $$ [\delta,x] \in {_sE} \quad\text{and}\quad \beta_s([\delta,x]) = [s,x] \in E_s. $$
            On the other hand, if $x \in X_t$, for some $t \in S^s$, then there is $z \in X$ such that $x = \alpha_t(z)$. Therefore,
                $$ [t,z] \in {_sE} \quad\text{and}\quad [a,y] = [s,x] = [s,\alpha_t(z)] \overset{\eqref{rel-2.2}}{=} [st,z]  = \beta_s([t,z]) \in E_s. $$
            That is, $E \setminus \iota(X) \subseteq E_1$. Hence, $E_0 = E_0 \cap \iota(X) \subseteq \iota(X_0) \subseteq E_0$. Since a partial action $(\alpha,X)$ is non-degenerate if and only if $X_0 = \emptyset$, we conclude that that $(\alpha,X)$ is non-degenerate if and only if $\emptyset = \iota(X_0) = E_0$ if and only if $(\beta,E)$ is non-degenerate.
        \end{proof}
    \end{prop}


\section{Application in semigroups and categories}

    \subsection{Partial semigroup actions.} In this subsection, we compare our definition of partial action when restricted to a semigroup with the one from \cite{kudryavtseva2023}. Furthermore, we show that the universal globalization constructed in Section \ref{sec:globalization} coincides with the tensor product globalization when the partial semigroup action is non-degenerate. 
    
    Throughout this subsection, $\mathcal{S}$ will denote a semigroup (that is, a semigroupoid such that $\mathcal{S}^{(2)} = \mathcal{S} \times \mathcal{S}$).

    \begin{defi} \cite[Definition 2.3]{kudryavtseva2023} \label{defi:semigroup_action}
        A \textbf{strong (left) partial semigroup action} of $\mathcal{S}$ on $X$ is a pair $(\{{_sX}\}_{s \in \mathcal{S}},\{\alpha_s\}_{s \in \mathcal{S}})$ consisting of a collection $\{_sX\}_{s \in \mathcal{S}}$ of subsets of $X$ and a collection $\{\alpha_s\}_{s \in \mathcal{S}}$ of functions $\alpha_s \colon {_sX} \to X$ such that:
        \begin{itemize}
            \item[(S)] $\alpha_t^{-1}(_sX \cap X_t) = {_tX} \cap {_{st}X}$ and, for $x \in {_tX} \cap {_{st}X}$, $(\alpha_s \circ \alpha_t)(x) = \alpha_{st}(x)$.
        \end{itemize}
        Such action is called \textbf{global} if ${_sX} = X$, for every $s \in \mathcal{S}$, and called \textbf{unitary} if $X = \cup_{s \in \mathcal{S}} X_s$.
    \end{defi}

    The strong partial actions described in \cite{kudryavtseva2023} are the type of partial actions that eventually can be globalizable \cite[Proposition 2.9]{kudryavtseva2023}. In \cite[Lemma 2.4]{kudryavtseva2023} it has been shown that strong partial monoid actions (in the sense of \cite[Definition 2.4]{hollings2007}) are precisely the strong and unitary partial semigroup actions of a monoid. We also observe that global semigroup actions are always non-degenerate, whereas partial actions are not necessarily so.

    \begin{prop} \label{prop:semigroup_actions}
        With the above notation, it follows that: \begin{itemize} \item[(i)] $\alpha$ is a strong partial semigroup action of $\mathcal{S}$ on $X$ (in the sense of Kudryavtseva and Laan) if and only if $\alpha$ is a partial semigroupoid action of $\mathcal{S}$ on $X$.
        
        \item[(ii)] $\alpha$ is a strong global semigroup action  of $\mathcal{S}$ on $X$ (in the sense of Kudryavtseva and Laan) if and only if $\alpha$ is a non-degenerate global semigroupoid action of $\mathcal{S}$ on $X$. \end{itemize}

        \begin{proof}
            Since $\mathcal{S}^{(2)} = \mathcal{S} \times \mathcal{S}$, condition (S) is equivalent to \eqref{pact-1} and \eqref{pact-2}. This proves (i).
            
            For (ii), assume that $\alpha$ is a strong global semigroup action of $\mathcal{S}$ on $X$. Then ${_sX} = X$, for all $s \in \mathcal{S}$. Thus, the action is non-degenerate, and $X={_sX} = {_tX} = {_{st}X}$. Thus, conditions \eqref{gact-1} and \eqref{gact-2} are trivially satisfied. Now assume that $\alpha$ is a non-degenerate global semigroupoid action of $\mathcal{S}$ on $X$. Since $\mathcal{S}^s = \mathcal{S}$, we obtain ${_sX} = {_tX}$, for all $s,t \in \mathcal{S}$. Define $X_1 = {_sX}$. From Lemma \ref{lema:global_charact}(ii) we obtain that $X_t \subseteq {_sX}$, for every $s,t \in \mathcal{S}$. Therefore
                $$ \cup_{s \in \mathcal{S}} ({_sX} \cup X_s) \subseteq \cup_{s \in \mathcal{S}} {_sX} = X_1. $$
            Since the action is non-degenerate, we have $X_1 = X$, which concludes the proof.
        \end{proof}
    \end{prop}

    Now, we relate our definition of (universal) globalization to the one presented for semigroups.

    Our notion of globalization is analogous to the respective notion from \cite[Definition 2.8]{kudryavtseva2023}. That is, a globalization of a partial semigroup action $(\alpha,X)$ is a triple $(i,\gamma,Y)$ where $(\gamma,Y)$ is a global semigroup action and $i \colon (\alpha,X) \to (\gamma,Y)$ is an embedding. On the other hand, the term \textit{universal globalization} was not adopted in \cite{kudryavtseva2023}, as two different globalizations of a partial semigroup action were considered, each one being \textit{universal} in a certain sense. The globalization we are interested in is called \textbf{the tensor product globalization}. It is constructed as follows.

    We start by extending the semigroup $\mathcal{S}$ to a monoid $\mathcal{S}^1$, by attaching a new element $1$ as a symbol outside of $\mathcal{S}$. The partial action of $\mathcal{S}$ on $X$ is then extended to an action of $\mathcal{S}^1$ on $X$ by defining ${_1X} = X$ and $\alpha_1 = id_X$. This extended action is both strong and unitary. Next, we define $\approx$ as the equivalence relation on $\mathcal{S}^1 \times X$ generated by the relation
    \begin{align*}
        (s,x) \sim (t,y) \text{ if there is } u \in \mathcal{S}^1 \text{ such that } s = ut, x \in {_uX} \text{ and } \alpha_u(x) = y. \tag{$\ast$} \label{rel:semigroup}
    \end{align*}
    Let $\mathcal{S}^1 \otimes X$ be the quotient set $(\mathcal{S}^1 \times X)/\approx$ and define $\gamma_s(t \otimes x) = st \otimes x$, for every $t \otimes x \in \mathcal{S}^1 \otimes X$. This is a global action of $\mathcal{S}$ on $\mathcal{S}^1 \otimes X$ and the function $i \colon X \to \mathcal{S}^1 \otimes X$, defined as $i(x) = 1 \otimes x$, is an embedding. This triple $(i,\gamma,\mathcal{S}^1 \otimes X)$ is the tensor product globalization of $(\alpha,X)$.
    
    It has been shown in \cite[Theorem 3.5]{kudryavtseva2023} that $(i,\gamma,\mathcal{S}^1 \otimes X)$ is an initial object in the \textit{category of globalizations} of $(\alpha,X)$, denoted by $\mathcal{G}(\alpha,X)$. The category $\mathcal{G}(\alpha,X)$ has as objects the embeddings $\iota \colon (\alpha,X) \to (\beta,Z)$, where $(\beta,Z)$ is a global semigroup action, and the morphisms between the objects $\iota \colon (\alpha,X) \to (\beta,Z)$ and $\varphi \colon (\alpha,X) \to (\theta,W)$ are the morphisms $\Phi \colon (\beta,Z) \to (\theta,W)$ such that $\Phi \circ \iota = \varphi$.

    \begin{obs} \label{obs:remark_semigroups}
        Let $(\alpha,X)$ be a partial semigroup action of $\mathcal{S}$. 
        \begin{itemize}
            \item[(1)] Let $\mathcal{A}^\ast(\mathcal{S})$ be the full subcategory of $\mathcal{A}(\mathcal{S})$ whose objects are the non-degenerated global actions of $\mathcal{S}$, and $I^\ast \colon \mathcal{A}^\ast(\mathcal{S}) \to \mathcal{A}_p(\mathcal{S})$ be the inclusion functor. Then $((\alpha,X) \downarrow I^\ast)$ is a proper full subcategory of $((\alpha,X) \downarrow I)$ (recall the definition in Proposition \ref{prop:universal_objeto_inicial}).

            \item[(2)] The category of globalizations $\mathcal{G}(\alpha,X)$ is a proper full subcategory of $((\alpha,X) \downarrow I^\ast)$. In the category $((\alpha, X) \downarrow I^\ast)$ the objects consist of all morphisms of partial actions $(\alpha, X) \to (\beta, Y)$, whereas $\mathcal{G}(\alpha, X)$ restricts its objects to embeddings only.
            
            \item[(3)] The proof given in \cite[Theorem 3.5]{kudryavtseva2023} does not use the hypothesis that an object $\iota \in \mathcal{G}(\alpha,X)$ is an embedding. Thus, the tensor product globalization $(i,\gamma, \mathcal{S}^1 \otimes X)$ also provides an initial object for the category $((\alpha,X) \downarrow I^\ast)$.
        \end{itemize}
    \end{obs}

    Due to Proposition \ref{prop:degenerate-1}, the universal globalization $(\iota,\beta,E)$ of a partial action $(\alpha,X)$, given by Theorem \ref{teo:globalization}, is non-degenerate if and only if $(\alpha,X)$ is non-degenerate. On the other hand, any global semigroup action (in the sense of Kudryavtseva and Laan) is non-degenerate. In particular, the tensor product globalization $(i,\gamma,\mathcal{S}^1 \otimes X)$ is non-degenerate. In the next result, we relate the globalizations $(\iota,\beta,E)$ and $(i,\gamma,\mathcal{S}^1 \otimes X)$.

    \begin{prop}
        Let $(\alpha,X)$ be a partial semigroupoid action of the semigroup $\mathcal{S}$. Then the tensor product globalization $(\beta,\mathcal{S}^1 \otimes X)$ and the universal globalization $(\beta,E)$ are isomorphic in $\mathcal{A}(\mathcal{S})$ if and only if $(\alpha,X)$ is a non-degenerate partial action.

        \begin{proof}
            Suppose that $(\alpha,X)$ is non-degenerate. From Proposition \ref{prop:degenerate-1}, it follows that $(\iota,\beta,Y)$ is non-degenerate and, hence, $(\iota,\beta,Y) \in ((\alpha,X) \downarrow I^\ast)$. By Proposition \ref{prop:universal_objeto_inicial}(i), $\iota$ is an initial object in the category $((\alpha,X) \downarrow I)$. Since $((\alpha,X) \downarrow I^\ast)$ is a full subcategory of $((\alpha,X) \downarrow I)$, it follows that $\iota$ is an initial object in the category $((\alpha,X) \downarrow I^\ast)$.

            On the other hand, $i \colon (\alpha,X) \to (\gamma,\mathcal{S}^1 \otimes X)$ is an initial object in $(\alpha(X) \downarrow I^\ast)$. Since initial objects are unique, up to isomorphism, there is an isomorphism $\Phi \colon \iota \to i$. By definition, $\Phi$ is an isomorphism of partial actions $(\beta,E) \to (\gamma,\mathcal{S}^1 \otimes X)$.

            To prove the converse, we will show that an isomorphism of partial actions preserve the non-degenerate part of the actions. In fact, let $(\alpha,X)$ and $(\beta,Y)$ be partial semigroup actions and suppose that $\varphi \colon (\alpha,X) \to (\beta,Y)$ is an isomorphism of partial actions. Then $\varphi^{-1} \colon (\beta,Y) \to (\alpha,X)$ is an isomorphism of partial actions. Hence, for each $s \in \mathcal{S}$, we have
            \begin{gather*}
                {_sY} = (\varphi \circ \varphi^{-1})({_sY}) \subseteq \varphi({_sX}) \subseteq {_sY}.
            \end{gather*}
            That is, $\varphi({_sX}) = {_sY}$. Hence,
            \begin{gather*}
                \varphi(X_s) = (\varphi \circ \alpha_s)({_sX}) = (\beta_s \circ \varphi)({_sX}) = \beta_s({_sY}) = Y_s.
            \end{gather*}
            Therefore,
                $$ \varphi(X_1) = \varphi \left( \cup_{s \in S} ({_sX} \cup X_s) \right) = \cup_{s \in S} (\varphi({_sX}) \cup \varphi(X_s)) = \cup_{s \in S} ({_sY} \cup Y_s) = Y_1. $$
            Now, suppose that $\varphi \colon (\beta,E) \to (\gamma,\mathcal{S}^1 \otimes X)$ is an isomorphism of global actions. Then
                $$ E = \varphi^{-1}( \mathcal{S}^1 \otimes X) = \varphi^{-1}( (\mathcal{S}^1 \otimes X)_1) = E_1. $$
            That is, $E_0 = E \setminus E_1 = \emptyset$. Therefore $(\beta,E)$ is non-degenerate and, from Proposition \ref{prop:degenerate-1}, it follows that $(\alpha,X)$ is non-degenerate.
        \end{proof}
    \end{prop}

    The difference between the universal globalization and the tensor product globalization lies in how we deal with the degenerate part of $(\alpha,X)$. Notice that, when $\mathcal{S}$ is a semigroup, the relation $R$ defined in Lemma \ref{lema:right_order} leads to $sRt$, for all $s,t \in \mathcal{S}$. That is, $\mathcal{S}/R$ has only one equivalence class, which we will denote by $\rho$. Therefore ${_{\rho}X} = \cup_{s \in \mathcal{S}} ({_sX} \cup X_s) = X_1$, and
        $$ D = (\{\delta\} \times X) \cup (\mathcal{S}\times X_1) = (\mathcal{S}^\delta \times X_1) \cup (\{\delta\} \times X_0). $$
    Consequently, for each $s \in \mathcal{S}$, we have ${_sE} = E \setminus \iota(X_0)$. On the other hand, since $(i,\gamma,\mathcal{S}^1 \otimes X)$ is a globalization of $(\alpha,X)$, there is a unique morphism of partial actions $\varphi \colon (\beta,E) \to (\gamma,\mathcal{S}^1 \otimes X)$ such that $\varphi \circ \iota = i$. This morphism is defined by
        $$ \varphi([\delta,x]) = i(x) = 1 \otimes x \quad\text{and}\quad \varphi([s,x]) = \gamma_s(i(x)) = s \otimes x. $$
    Therefore, $\varphi \colon E \to \mathcal{S}^1 \otimes X$ is a bijection. However, the inverse function $\varphi^{-1} \colon \mathcal{S}^1 \otimes X \to E$ may not be a morphism of partial actions, since in $(\gamma,\mathcal{S}^1 \otimes X)$ every element acts on $i(X_0) = \varphi(\iota(X_0))$ as the identity function, while in $(\beta,E)$ no element can act on $\iota(X_0)$.\\

    We end this section noticing that the fact that any degenerate partial semigroup action admits a non-degenerate globalization that provides an initial object in the category $((\alpha,X) \downarrow I^\ast)$ is a convenient property of semigroups. Indeed, if $(\beta,Y)$ is a non-degenerate global action of $\mathcal{S}$, then ${_sY} = Y$ for every $s \in \mathcal{S}$. Therefore, in order to define a morphism $\Phi \colon (\beta,E) \to (\theta,Y)$, there is no need to verify that $\Phi({_sE}) \subseteq {_sY}$. The next example illustrates, in the context of semigroupoids, how a degenerate partial action may fail to have such a non-degenerate globalization.

    \begin{exe}
        Let $S = \{e,f\}$, $S^{(2)} = \{(e,e),(f,f)\}$ and $ee = e$, $ff = f$. Then $S$ is a semigroupoid (furthermore, $S$ is a groupoid, but this will not be relevant). Let $Y = \{0,1,2,3\}$ and define
        \begin{align*}
            {_eY} &= \{0,1,2\}, & \beta_e(0) = 0,\ \beta_e(1) = 1,\ \beta_e(2) = 0, \\
            {_fY} &= \{0,1,3\}, & \beta_f(0) = 0,\ \beta_f(1) = 1,\ \beta_f(3) = 0.
        \end{align*}
        Then $(\beta,Y)$ is a non-degenerate global action of $G$. Let $X = \{1,2\}$ and $\alpha$ be the restriction of $\beta$ to $X$, as in Proposition \ref{prop:restriction}. Then
        \begin{gather*}
            {_eX} = \{1\}, \quad \beta_e(1) = (1), \quad {_fX} = \{1\} \quad\text{and}\quad \beta_f(1) = (1).
        \end{gather*}
        In this case $(\alpha,X)$ is a degenerate global action of $S$. We know that $i \colon X \to Y$, given by $i(1) = 1$ and $i(2) = 2$ is an embedding. Notice that $j \colon X \to Y$, given by $j(1) = 1$ and $j(2) = 3$, is also an embedding.
        
        Suppose that $(\iota,\theta,E)$ is a non-degenerate globalization of $(\alpha,X)$. Since $S^s \neq \emptyset$ for every $s \in S$, we obtain from Lemma \ref{lema:global_charact}(ii) that $E_s \subseteq {_tE}$, for some $t \in S^s$. Therefore, $E \subseteq \cup_{s \in S} {_sE}$. Now, the embedding $\iota \colon X \to E$ must satisfy $\iota(2) \in {_eE}$ or $\iota(2) \in {_fE}$. Suppose further that $(\iota,\theta,E)$ is initial object in the category $((\alpha,X) \downarrow I^\ast)$. Then there must be two morphisms of partial actions, $i', j' \colon (\theta,E) \to (\beta,Y)$, satisfying $i' \circ \iota = i$ and $j' \circ \iota = j$. Therefore,
            $$ i'(\iota(2)) = i(2) = 2 \in {_eY} \setminus {_fY} \quad\text{and}\quad j'(\iota(2)) = j(2) = 3 \in {_fY} \setminus {_eY}. $$
        Since $i'$ is a morphism of partial actions, if $\iota(2) \in {_fE}$ then $i'(\iota(2)) \in {_fY}$, a contradiction. Hence, it must be $\iota(2) \in {_eE}$. But in this case, as $j'$ is a morphism of partial actions, it must be $j'(\iota(2)) \in {_eY}$, a contradiction. Therefore, the degenerate partial action $(\alpha,X)$ does not admit a non-degenerate globalization that provides an initial object for $((\alpha,X) \downarrow I^\ast)$.
    \end{exe}

    \subsection{Partial category actions.} In this subsection, we investigate the relation between partial category actions and partial semigroupoid actions of a category and show that our construction of universal globalization coincides with the one in \cite[\textsection 3]{nystedt2018}.
    
    Throughout this subsection, $\mathcal{C}$ will denote a category (recall our conventions in Definition \ref{defi:category}). The following characterization of partial category actions was given in \cite[Proposition 8]{nystedt2018}.

    \begin{defi}
        A \textbf{partial action} of $\mathcal{C}$ on a set $X$ is a pair $(\{{_gX}\}_{g \in \mathcal{C}}, \{\alpha_g\}_{g \in \mathcal{C}})$ consisting of a collection $\{{_gX}\}_{g \in \mathcal{C}}$ of subsets of $X$ and a collection $\{\alpha_g\}_{g \in \mathcal{C}}$ of functions $\alpha_g \colon {_gX} \to X$ such that:
        \begin{enumerate} \Not{c}
            \item $X = \cup_{e \in \mathcal{C}_0} {_eX}$ and $\alpha_e = id_{_eX}$, for all $e \in \mathcal{C}_0$; \label{cat_action-1}
            \item If $g \in \mathcal{C}$, then ${_gX} \subseteq {_{D(g)}X}$; \label{cat_action-2}
            \item If $gh$ exists, then $\alpha_h^{-1}(X_h \cap {_gX}) = {_hX} \cap {_{gh}X}$ and, for $x \in {_hX} \cap {_{gh}X}$, $(\alpha_g \circ \alpha_h)(x) = \alpha_{gh}(x)$.\label{cat_action-3}
        \end{enumerate}
        Such action is called \textbf{global} if, furthermore, it satisfies:
        \begin{enumerate} \Not{c} \setcounter{enumi}{3}
            \item ${_gX} = {_{D(g)}X}$, for every $g \in \mathcal{C}$. \label{cat_action-4}
        \end{enumerate}
    \end{defi}

    Compared to the original characterization, we omit the collection of subsets $\{X_g\}_{g \in \mathcal{C}}$, since it is implicitly defined as the images of the collection of functions $\{\alpha_g\}_{g \in \mathcal{C}}$, and interchange the roles of $X_g$ and ${_gX}$. This change aims the comparison with partial group actions \cite[Definition 1.2]{exel1998}, where $\alpha_g \colon X_{g^{-1}} \to X_g$, or in our notation, $X_g = X_g$ and ${_gX} = X_{g^{-1}}$. It is straightforward from \eqref{cat_action-1} that any partial category action is non-degenerate.\\

    Similarly to the relation between semigroups and monoids (\cite[Lemma 2.4]{kudryavtseva2023}), there is a particular class of partial semigroupoid actions that correspond to partial category actions. To describe it, we need the notion of an \textit{identity element} in a semigroupoid.

    \begin{defi}
        Let $S$ be a semigroupoid. An element $e \in S$ is called an \textbf{identity} if $(e,e) \in S^{(2)}$ and $se = s$ and $et = t$, whenever $(s,e),(e,t) \in S^{(2)}$. We denote by $S_0$ the set of identities of $S$.
    \end{defi}

    Notice that $S_0$ may be the empty set. Furthermore, for each $s \in S$, there exists at most one identity $e \in S^s$. Indeed, if $e,f \in S_0 \cap S^s$, then $s = se = (sf)e$. Hence, $fe$ is defined, and since both $e$ and $f$ are identities, we conclude that $e = fe = f$.

    \begin{lema} \label{lema:identity_action}
        Let $S$ be a semigroupoid, $(\alpha,X)$ be a partial semigroupoid action, and $e \in S_0$. It holds that:
        \begin{itemize}
            \item[(i)] $\alpha_e$ is a projection; that is, $X_e \subseteq {_eX}$ and $\alpha_e(x) = x$, for every $x \in X_e$.
            \item[(ii)] If $(e,s) \in S^{(2)}$, then $X_s \subseteq X_e$.
        \end{itemize}

        \begin{proof}
            Let us start with (ii). Suppose that $(e,s) \in S^{(2)}$. Then $es = s$, and by \eqref{pact-1} we obtain $\alpha_s^{-1}({_eX} \cap X_s) = {_sX} \cap {_{es}X} = {_sX}$. Applying $\alpha_s$ to this equality, we obtain ${_eX} \cap X_s = X_s$. Therefore $X_s \subseteq {_eX}$. Let $x \in {_sX} = {_sX} \cap {_{es}X}$. Then condition \eqref{pact-2} implies $(\alpha_e \circ \alpha_s)(x) = \alpha_{es}(x) = \alpha_s(x)$. That is, for every $y = \alpha_s(x) \in X_s$ we have $\alpha_e(y) = y \in X_e$. Thus, $X_s \subseteq {_eX} \cap X_e$.

            Applying the previous argument for $s=e$ we obtain $X_e \subseteq {_eX} \cap X_e$ and $(\alpha_e \circ \alpha_e)(x) = \alpha_e(x)$, for every $x \in {_eX}$. Therefore, $X_e \subseteq {_eX}$, from which we obtain $X_s \subseteq {_eX} \cap X_e = {_eX}$. This concludes (ii) and proves (i).
        \end{proof}
    \end{lema}

    When $S = \mathcal{C}$ is a category, Lemma \ref{lema:identity_action} implies that $X_g \subseteq {_{R(g)}X}$, for every $g \in \mathcal{C}$, and for each $e \in \mathcal{C}_0$ the function $\alpha_e$ is a projection. The next example shows that, for a partial semigroupoid action of a category, conditions \eqref{cat_action-1} and \eqref{cat_action-2} may not hold.

    \begin{exe}
        Let $\mathcal{C}$ be the category $\{R(g),g,D(g)\}$ consisting of two objects and one (non identity) morphism between them. Let $X = \{0,1,2\}$ and consider the subsets
        \begin{align*}
            {_gX} = \{0,2\}, \quad {_{D(g)}X} = \{0\} \quad\text{and}\quad {_{R(g)}X} = \{0,1\},
        \end{align*}
        and let $\alpha_h \equiv 0$, for every $h \in \mathcal{C}$. Then $(\alpha,X)$ is a partial semigroupoid action of $\mathcal{C}$. However, $(\alpha,X)$ is not a partial category action of $\mathcal{C}$ since $\alpha_{R(g)} \neq id_{_{R(g)}X}$ and ${_gX} \not\subseteq {_{D(g)}X}$.
    \end{exe}

    \begin{defi}
        Let $S$ be a semigroupoid and $(\alpha,X)$ be a partial semigroupoid action. We say that $(\alpha,X)$ is a \textbf{categorical} partial action of $S$ if the following conditions hold:
        \begin{enumerate} \Not{C}
            \item $X = \cup_{s \in S}\ X_s$. \label{categorical-1}
            \item $\alpha_e = id_{_eX}$, for every $e \in S_0$. \label{categorical-2}
            \item If $e \in S_0$ and $(s,e) \in S^{(2)}$, then ${_sX} \subseteq {_eX}$. \label{categorical-3}
        \end{enumerate}
    \end{defi}

    Observe that in \eqref{categorical-1}, we ask for $X = \cup_{s \in S} X_s$ instead of $X = \cup_{e \in S_0} {_eX}$. By doing so, the categorical partial actions generalize the unitary partial actions of a semigroup (Definition \ref{defi:semigroup_action}).

    \begin{prop} \label{prop:act-cat-sgpd}
        Let $\mathcal{C}$ be a category, $X$ be a set, $\{{_gX}\}_{g \in \mathcal{C}}$ be a family of subsets of $X$ and $\{\alpha_g\}_{g \in \mathcal{C}}$ be a family of functions $\alpha_g \colon {_gX} \to X$. Then the following are equivalent:
        \begin{itemize}
            \item[(i)] $(\alpha,X)$ is a partial (global) category action of $\mathcal{C}$ (in the sense of Nystedt).
            \item[(ii)] $(\alpha,X)$ is a categorical partial (global) semigroupoid action of $\mathcal{C}$.
        \end{itemize}

        \begin{proof}
            We will first prove the equivalence for partial actions. Notice that \eqref{pact-1} and \eqref{pact-2} are equivalent to \eqref{cat_action-3}; \eqref{categorical-2} is equivalent to the second part of \eqref{cat_action-1}; and \eqref{categorical-3} is equivalent to \eqref{cat_action-2}.

            Assume (i). To prove (ii) it remains to show \eqref{categorical-1}. From the second part of \eqref{cat_action-1} we have $\alpha_e = id_{_eX}$, for every $e \in \mathcal{C}_0$. That is, ${_eX} = X_e$. From the first part of \eqref{cat_action-1} we obtain
                $$ X = \cup_{e \in \mathcal{C}_0}\ {_eX} = \cup_{e \in \mathcal{C}_0} X_e \subseteq \cup_{g \in \mathcal{C}} X_g \subseteq X. $$
            Therefore $X = \cup_{g \in \mathcal{C}} X_g$. Conversely, assume (ii). It remains to show the first part of \eqref{cat_action-1}, that is, $X = \cup_{e \in \mathcal{C}_0} {_eX}$. From Lemma \ref{lema:identity_action}(ii) and \eqref{categorical-2} we have $X_g \subseteq {X_{R(g)}} = {_{R(g)}X}$, for every $g \in C$. From \eqref{categorical-1} we obtain
                $$ X = \cup_{g \in \mathcal{C}} X_g \subseteq \cup_{g \in \mathcal{C}} {_{R(g)}X} = \cup_{e \in \mathcal{C}_0} {_eX} \subseteq X. $$
            Therefore $X = \cup_{e \in \mathcal{C}_0} {_eX}$.

            For the equivalence on global actions it remains to prove that \eqref{gact-1} and \eqref{gact-2} are equivalent to \eqref{cat_action-4}. Moreover, we will show that \eqref{gact-2} implies \eqref{cat_action-4}, which implies \eqref{gact-1}. Since $D(g) \in \mathcal{C}^g$, for every $g \in \mathcal{C}$, we obtain from Lemma \ref{lema:G1implicaG2} that \eqref{gact-1} implies \eqref{gact-2}. Therefore, over a category all three conditions are equivalent.

            Suppose that \eqref{gact-2} holds, that is, ${_hX} = {_{gh}X}$, whenever the composition $gh$ is defined. In particular, for $h = D(g)$, the composition $gD(g)$ exists and, hence, ${_{D(g)}X} = {_{gD(g)}X} = {_gX}$. This proves \eqref{cat_action-4}. Assume \eqref{cat_action-4} and notice that, for $g \in \mathcal{C}$,
                $$ \mathcal{C}^g = \{h \in \mathcal{C} \colon \exists gh\} = \{h \in \mathcal{C} \colon R(h) = D(g)\}. $$
            Therefore, if $\mathcal{C}^g = \mathcal{C}^h$, then $D(g) = D(h)$. In this case ${_gX} = {_{D(g)}X} = {_hX}$. This proves that \eqref{gact-1} is satisfied.
        \end{proof}
    \end{prop}

    Now we relate our notion of (universal) globalization to the one presented for categories.

    The notion of globalization used in \cite[Definition 11]{nystedt2018} is more general than ours. The term \textit{globalization} means any morphism of partial actions $\iota \colon (\alpha,X) \to (\beta,Y)$ where $(\beta,Y)$ is a global action, while we require $\iota$ to be an embedding. On the other hand, \textit{universal globalization} means a globalization $(\iota,\beta,Y)$ such that, for every morphism $\varphi \colon (\alpha,X) \to (\theta,Z)$ where $(\theta,Z)$ is a global action, there is a unique morphism $\Phi \colon (\beta,Y) \to (\theta,Z)$ satisfying $\Phi \circ \iota = \phi$. Apart from the definition of globalization, the idea of \textit{universal globalization} agree with ours.

    \begin{obs} \label{obs:cat}
        Let $(\alpha,X)$ be a categorical partial action of the category $\mathcal{C}$.
        \begin{itemize}
            \item[(1)] In \cite[Theorem 4]{nystedt2018} it has been proved that $(\alpha,X)$ admits a universal globalization $(i,\gamma,Y)$. Furthermore, \cite[Remark 22]{nystedt2018} shows that $i$ is, in fact, an embedding.
            
            \item[(2)] From (1) and Proposition \ref{prop:universal_objeto_inicial}(ii), it follows that, for any universal globalization $(\iota,\beta,Z)$ of $(\alpha,X)$, $\iota$ must be an embedding. Therefore, a universal globalization of $(\alpha,X)$ is a globalization, in the sense of Definition \ref{def:globalization}.
            
            \item[(3)] Let $C\mathcal{A}(\mathcal{C})$ be the full subcategory of $\mathcal{A}_p(\mathcal{C})$ whose objects are the categorical global actions of $\mathcal{C}$, and $I_C \colon C\mathcal{A}(\mathcal{C}) \to \mathcal{A}_p(\mathcal{C})$ be the inclusion functor. Then the category $((\alpha,X) \downarrow I_C)$ is a full subcategory of the category $((\alpha,X) \downarrow I)$ (recall the definition in Proposition \ref{prop:universal_objeto_inicial}).

            \item[(4)] From \cite[Theorem 4]{nystedt2018}, (2) and Proposition \ref{prop:universal_objeto_inicial}(i), it follows that the universal globalization $(i,\gamma,Y)$ is an initial object in the category $((\alpha,X) \downarrow I_C)$.
        \end{itemize}
    \end{obs}

    The universal globalization from \cite{nystedt2018} is constructed as follows. Given a partial category action $(\alpha,X)$ of $\mathcal{C}$, let $\overline{X} = \{(g,x) \colon x \in {_{D(g)}X}\}$. For $(g,x),(g',x') \in \overline{X}$, denote $(g,x) - (g',x')$ if and only if one of the following conditions is satisfied:
    \begin{enumerate} \Not{$\ast$}
        \item There is $h \in \mathcal{C}$ such that $g'h$ is defined, $g'h = g$, $x \in {_hX}$ and $\alpha_h(x) = x'$. \label{rel-cat-1}
        \item $x = x'$ and $g,g' \in \mathcal{C}_0$ (and in this case $x \in {_gX} \cap {_{g'}X}$). \label{rel-cat-2}
    \end{enumerate}
    Let $\simeq$ be the smallest equivalence relation on $\overline{X}$ that contains the relation $-$ and denote $Y = X/\simeq$. For each $g \in \mathcal{C}$, let $(h,x) \in {_g\overline{X}}$ if and only if $(h,x) \simeq (h',x')$ for some $h' \in \mathcal{C}$ such that $gh'$ is defined. Denote ${_gY} = {_g\overline{X}}/\simeq$ and define $\gamma_g \colon {_gY} \to Y$ by $\gamma_g([h,x]) = [gh',x']$. Then $(\gamma,Y)$ is a global action. From condition \eqref{cat_action-1}, for every $x \in X$ there is $e \in \mathcal{C}_0$ such that $x \in {_eX}$. Define $i \colon X \to Y$ by $i(x) = [e,x]$. The triple $(i,\gamma,Y)$ is the desired universal globalization for $(\alpha,X)$.\\

    To prove that our universal globalization $(\iota,\beta,E)$ and the one described above are isomorphic, we will use Remark \ref{obs:cat}(4) and the uniqueness, up to isomorphism, of initial objects.

    \begin{prop} \label{prop:global-categorical}
        Let $(\alpha,X)$ be a categorical partial action of $\mathcal{C}$. Then the universal globalization $(\iota,\beta,E)$ of Theorem \ref{teo:globalization} is a categorical action of $\mathcal{C}$.

        \begin{proof}
            We start by summarizing the construction of $(\iota,\beta,E,)$ when $\mathcal{C}$ is a category and $(\alpha,X)$ is a categorical partial action of $\mathcal{C}$. The relation $R$ defined in Lemma \ref{lema:right_order} leads to $gRh$ if and only if $D(g) = D(h)$. Denote by $\rho_g = \rho_{D(g)}$ the equivalence class of $g$. Then, for each $e \in \mathcal{C}_0$, we have
                $$ {_{\rho_e}X} = \bigcup \left( \{{_hX} \colon D(h) = e\} \cup \{X_h \colon R(h) = e\} \right). $$
            Notice that, if $R(h) = e$, then $eh$ is defined. Therefore, it follows from Lemma \ref{lema:identity_action}(ii) and \eqref{categorical-2} that $X_h \subseteq X_e = {_eX}$. On the other hand, if $D(h) = e$ , then \eqref{categorical-3} implies ${_hX} \subseteq {_eX}$. Therefore ${_{\rho_e}X} \subseteq {_eX}$. The reverse inclusion is trivial, since $D(e) = e$. Consequently, the set $D$ is defined as
                $$ D = (\{\delta\} \times X) \cup \{(g,x) \colon x \in {_{D(g)}X}\}. $$
            Let $E = D/\approx$, where $\approx$ is the equivalence relation given in Definition \ref{defi:rel-2}. For an element $(a,x) \in D$ we will denote its equivalence class by $[a,x] \in E$. Now we prove conditions \eqref{categorical-1}, \eqref{categorical-2} and \eqref{categorical-3}.

            First, let $(g,x) \in D$ with $g \in \mathcal{C}$. Since the composition $R(g)g$ is defined, we have $[g,x] \in {_{R(g)}E}$ by definition. Now, let $(\delta,x) \in D$. Since $(\alpha,X)$ is a categorical partial action, it follows from \eqref{categorical-1} that $x \in X_g$ for some $g \in \mathcal{C}$. Since $R(g)g = g$, Lemma \ref{lema:identity_action}(ii) implies $X_g \subseteq X_{R(g)}$ and, from \eqref{categorical-2}, we obtain $X_{R(g)} = {_{R(g)}X}$. Therefore $(R(g),x) \in D$ and $(\delta,x) \sim (R(g),x)$ by \eqref{rel-2.1}. From the first case, we have $[\delta,x] = [R(g),x] \in {_{R(g)}E}$. Therefore $E = \cup_{e \in \mathcal{C}_0}\ {_eE}$.

            Let $[a,x] \in {_eE}$ for some $e \in \mathcal{C}_0$. Then either $[a,x] = [t,y]$ for some $t \in \mathcal{C}$ such that $R(t) = e$; or $a = \delta$ and $x \in {_eX}$. In the first case we have
                $$ \beta_e([a,x]) = [et,y] = [t,y] = [a,x], $$
            while in the second, the previous argument shows that
                $$ \beta_e([\delta,x]) = [e,x] = [\delta,x]. $$
            Therefore $\beta_e = id_{_eX}$. This proves \eqref{categorical-2}. Furthermore, we obtain that ${_eE} = E_e$ and, hence, $E = \cup_{e \in \mathcal{C}_0}\ {_eE} = \cup_{e \in \mathcal{C}_0} E_e$. This proves \eqref{categorical-1}. Condition \eqref{categorical-3} follows from $(\beta,E)$ being a global action. In fact, for each $g \in \mathcal{C}$, the unique $e \in \mathcal{C}_0$ such that $(g,e) \in \mathcal{C}^{(2)}$ is $e = D(g)$. Therefore,
                $$ {_gE} = {_{gD(g)}E} \overset{\eqref{gact-2}}{=} {_{D(g)}E}. $$
            This concludes that $(\beta,E)$ is a categorical action of $\mathcal{C}$.
        \end{proof}
    \end{prop}

    \begin{coro}
        Let $(\alpha,X)$ be a categorical partial action of $\mathcal{C}$, $(\iota,\beta,E)$ be the universal globalization of Theorem \ref{teo:globalization}, and $(i,\gamma,Y)$ be the universal globalization of \cite[Theorem 4]{nystedt2018}. Then $(\beta,E)$ and $(\gamma,Y)$ are isomorphic in the category $\mathcal{A}(\mathcal{C})$.

        \begin{proof}
            We know that $(\iota,\beta,E)$ is an initial object in the category $((\alpha,X) \downarrow I)$. Since $(\alpha,X)$ is a categorical partial action, it follows from Proposition \ref{prop:global-categorical} that $(\beta,E)$ is a categorical global action. Therefore $(\iota,\beta,E) \in C\mathcal{A}(\mathcal{C})$ and, hence, $(\iota,\beta,E)$ is an initial object in the category $((\alpha,X) \downarrow I_C)$.

            On the other hand, we obtain from Remark \ref{obs:cat} that $(i,\gamma,Y)$ is an initial object in the category $((\alpha,X) \downarrow I_C)$. Since initial objects are unique, up to isomorphism, there is an isomorphism $\varphi \colon \iota \to i$ in the category $((\alpha,X) \downarrow I_C)$. By definition, $\varphi$ is an isomorphism of partial actions $\varphi \colon (\beta,E) \to (\gamma,Y)$ such that $\varphi \circ \iota = i$. In particular $(\beta,E)$ is isomorphic to $(\gamma,Y)$ in $\mathcal{A}(\mathcal{C})$.
        \end{proof}
    \end{coro}


\section{A note on inverse semigroupoids} The globalization problem for partial actions of inverse semigroupoids on sets was studied in \cite{demeneghi}. However, the notions of partial and global actions for inverse semigroupoids differ slightly from our definitions of partial and global actions for semigroupoids. As a consequence, our construction of a universal globalization does not extend the results of \cite{demeneghi}. In this section, we present examples that highlight the differences between partial and global actions of inverse semigroupoids and semigroupoids.

\begin{defi}
    An \textbf{inverse semigroupoid} is a semigroupoid $S$ such that, for every $s \in S$, there is an unique element $s^{-1} \in S$ satisfying
        $$ (s,s^{-1}), (s^{-1},s) \in S^{(2)}, \quad ss^{-1}s = s \quad\text{and}\quad s^{-1}ss^{-1} = s. $$
    In this case, the element $s^{-1}$ is called the \textbf{inverse} of $s$.
\end{defi}

\begin{obs} \label{obs:inv-comute}
    Let $S$ be an inverse semigroupoid. Then, for all $s \in S$, we have
        $$ ss^{-1}, s^{-1}s \in E(S) = \{ e \in S \colon (e,e) \in S^{(2)} \text{ and } ee = e \}. $$
    It follows from \cite[Theorem 2.25]{cordeiro2019} that, if $e,f \in E(S)$ and $(e,f) \in S^{(2)}$, then $ef = fe$.
\end{obs}

An inverse semigroupoid in which every element has associated source and range identities is called an \textbf{inverse category}. The following example presents a simple inverse category that will be useful for illustrating the ideas developed in this subsection.

\begin{exe} \label{exe:0e1}
    Consider the semigroup $S = \{0, e, 1\}$ with multiplication defined by $0s = 0 = s0$ and $1s = s = s1$ for all $s \in S$, and $ee = e$. Then $S$ is an inverse semigroup, with $s^{-1} = s$ for every $s \in S$. Moreover, $S$ is a monoid with identity $1$, and therefore an inverse category.
\end{exe}

\begin{defi} \cite[Definition 2.4]{demeneghi}
    A \textbf{partial action} $\theta$ of an inverse semigroupoid $S$ on a set $X$ is a pair $\theta = (\{{X_s}\}_{s \in S},\{\theta_s\}_{s \in S})$ consisting of a collection $\{X_s\}_{s \in S}$ of subsets of $X$ and a collection $\{\theta_s\}_{s \in S}$ of maps $\theta_s \colon X_{s^{-1}} \to X_s$ such that
    \begin{enumerate} \Not{I}
        \item $\theta_e = id_{X_e}$ for all $e \in E(S)$. Moreover, for every $x \in X$, there exists $e \in E(S)$ such that $x \in X_e$. \label{I1}
        \item $X_s \subseteq X_{ss^{-1}}$, for every $s \in S$. \label{I2}
        \item $\theta_t^{-1}(X_t \cap X_{s^{-1}}) = X_{(st)^{-1}} \cap X_{t^{-1}}$, for all $(s,t) \in S^{(2)}$, and $\theta_s(\theta_t(x)) = \theta_{st}(x)$, for all $x \in X_{(st)^{-1}} \cap X_{t^{-1}}$. \label{I3}
    \end{enumerate}
    We say that $\theta$ is a \textbf{global action} if moreover it satisfies:
    \begin{enumerate} \Not{I} \setcounter{enumi}{3}
        \item $X_s = X_{ss^{-1}}$, for all $s \in S$. \label{I4}
    \end{enumerate}
\end{defi}

\begin{obs} \label{obs:inverse}
    Every partial inverse semigroupoid action $\theta = (\{{X_s}\}_{s \in S},\{\theta_s\}_{s \in S})$ is a partial semigroupoid action in which the domain of $\theta_s$ is ${_sX} = X_{s^{-1}}$ and the image of $\theta_s$ is $X_s = X_s$. With this notation, \eqref{I3} is equivalent to conditions \eqref{pact-1} and \eqref{pact-2}.
\end{obs}

The following two examples show that, for a general inverse category $S$, partial inverse semigroupoid actions of $S$ and partial category actions of $S$ do not necessarily coincide.

\begin{exe} \label{inverse-act-1}
    Let $S = \{0,e,1\}$ be the inverse monoid from Example \ref{exe:0e1}. Then every partial inverse semigroupoid action of $S$ on a set $X$ is, in fact, a global inverse semigroupoid action. Indeed, since $s^{-1} = s$ and $ss = s$ for every $s \in S$, the additional condition $X_s = X_{ss^{-1}}$ is automatically satisfied.
    
    In particular, let
        $$ X_0 = \{0\}, \quad X_e = \{0,e\}, \quad X_1 = \{0,e,1\} \quad\text{and}\quad \theta_s = id_{X_s}, \ \forall s \in S. $$
    Then $\theta = (\{X_s\}_{s \in S}, \{\theta_s\}_{s \in S})$ is a global inverse semigroupoid action of $S$ on $X=S$. 

    On the other hand, $\theta$ can also be viewed as a partial category action of $S$ on $X$, which is \emph{not} a global category action, since $X_0 \neq X_1 = X_{D(0)}$. That is, $\theta$ fails to satisfy condition \eqref{cat_action-4}.
\end{exe}

\begin{exe} \label{inverse-act-2}
    Let $S = \{0,e,1\}$ be the inverse monoid from Example \ref{exe:0e1}. Recall from Example \ref{exe:regular_action} that left multiplication on $S$ defines a global semigroupoid action of $S$ on $X = S$. In this setting, we have ${_0X} = S = {_eX} = {_1X}$, and the action is given by
    \begin{gather*}
        \alpha_0(0) = 0, \ \alpha_0(e) = 0, \ \alpha_0(1) = 0, \\
        \alpha_e(0) = 0, \ \alpha_e(e) = e, \ \alpha_e(1) = e, \\
        \alpha_1(0) = 0, \ \alpha_1(e) = e, \ \alpha_1(1) = 1.
    \end{gather*}
    Hence, $\alpha = (\{X_s\}_{s \in S}, \{\alpha_s\}_{s \in S})$ is a categorical global action of $S$ on $X$, and therefore a global category action of $S$. However, the same pair $\alpha$ cannot be regarded as a partial inverse semigroupoid action of $S$, since $0 \in E(S)$ but $\alpha_0 \neq id_{X_0}$. In other words, $\alpha$ fails to satisfy condition \eqref{I1}. This example thus illustrates that a global category action of an inverse category need not correspond to a partial inverse semigroupoid action.
\end{exe}

Recall from Proposition \ref{prop:act-cat-sgpd} that partial category actions of a category $\mathcal{C}$ are precisely those partial semigroupoid actions of $\mathcal{C}$ that are also categorical. In other words, for a given category $\mathcal{C}$, there exist more partial actions of $\mathcal{C}$ when it is regarded as a semigroupoid than when it is regarded as a category. Nevertheless, the partial actions of $\mathcal{C}$ as a category can be recovered by imposing additional requirements, namely conditions \eqref{categorical-1}, \eqref{categorical-2}, and \eqref{categorical-3}.

In contrast, it is not possible to formulate a set of conditions $\mathfrak{C}$ such that the partial inverse semigroupoid actions of an inverse semigroupoid $S$ correspond exactly to the partial semigroupoid actions of $S$ satisfying $\mathfrak{C}$. This is because condition \eqref{gact-1} is slightly more restrictive than condition \eqref{I4}, as will be shown in the next result.

\begin{prop} \label{prop:inverse}
   Let $S$ be an inverse semigroupoid, and let $\alpha$ be a partial inverse semigroupoid action of $S$ on a set $X$. Then the following conditions are equivalent:
    \begin{itemize}
        \item[(a)] $\alpha$ is a global semigroupoid action of $S$ on $X$.
        \item[(b)] $\alpha$ is a global inverse semigroupoid action of $S$ on $X$, and for every $e, f \in E(S)$ such that $(e, f) \in S^{(2)}$, we have $X_e = X_f$.
    \end{itemize}

    \begin{proof}
        Note first that $s^{-1} \in S^s$ for every $s \in S$. Hence, by Lemma \ref{lema:G1implicaG2}, condition \eqref{gact-1} implies condition \eqref{gact-2}. Therefore, it is enough to prove that $\alpha$ satisfies \eqref{gact-1} if and only if it satisfies (b). 
        
        From Remark \ref{obs:inv-comute} and the associativity of $S$, we know that for all $e, f \in E(S)$ with $(e, f) \in S^{(2)}$,
            $$ (f,s) \in S^{(2)} \iff (ef,s) \in S^{(2)} \iff (fe,s) \in S^{(2)} \iff (e,s) \in S^{(2)}. $$
        Thus, $S^e = S^f$. 
        
        Suppose that $\alpha$ satisfies \eqref{gact-1} and let $e, f \in E(S)$ be such that $(e, f) \in S^{(2)}$. Since $S^e = S^f \neq \emptyset$, we have
            $$ X_e \overset{\eqref{I1}}{=} {_eX} \overset{\eqref{gact-1}}{=} {_fX} \overset{\eqref{I1}}{=} X_f. $$
        Furthermore, since $s^{-1} = s^{-1}(ss^{-1})$, it follows that $S^{s^{-1}} = S^{ss^{-1}} \neq \emptyset$, and consequently
            $$ X_s = {_{s^{-1}}X} \overset{\eqref{gact-1}}{=} {_{ss^{-1}}X} = X_{ss^{-1}}. $$
        This proves (b). 
        
        Conversely, suppose that $\alpha$ satisfies (b). If $S^s = S^t \neq \emptyset$, then $t^{-1}t \in S^t = S^s = S^{s^{-1}s}$. Therefore, $(s^{-1}s,t^{-1}t) \in S^{(2)}$. Since $s^{-1}s, t^{-1}t \in E(S)$, it follows that
            $$ {_sX} = X_{s^{-1}} \overset{\eqref{I4}}{=} X_{s^{-1}s} \overset{(b)}{=} X_{t^{-1}t} \overset{\eqref{I4}}{=} X_{t^{-1}} = {_tX}. $$
        This proves that \eqref{gact-1} is satisfied. Hence, $\alpha$ is a global semigroupoid action.
    \end{proof}
\end{prop}

From the discussion above, we conclude that partial and global inverse semigroupoid actions, as defined in \cite{demeneghi}, cannot be obtained as particular cases of partial and global semigroupoid actions. Moreover, Examples \ref{inverse-act-1} and \ref{inverse-act-2} show that, for an inverse category $S$, the partial and global actions of $S$ as a category differ from the corresponding partial and global actions of $S$ as an inverse semigroupoid. Therefore, even if the notion of global semigroupoid action is modified, it does not appear possible to define a unified framework for partial and global actions of semigroupoids that simultaneously encompasses the cases of categories and inverse semigroupoids.\\

To conclude this section, we include an example illustrating the behavior of our universal globalization construction when applied to a partial inverse semigroupoid action. Let $S$ be an inverse semigroupoid and $(\alpha, X)$ a partial inverse semigroupoid action. By Remark \ref{obs:inverse}, $(\alpha, X)$ is, in particular, a partial semigroupoid action of $S$. Then, by Theorem \ref{teo:globalization}, we obtain a global semigroupoid action $(\beta, E)$. However, the action $(\beta, E)$ need not be an inverse semigroupoid action, and the embedding $\iota \colon (\alpha, X) \to (\beta, E)$ may fail to factor every morphism $(\alpha, X) \to (\gamma, Y)$, where $(\gamma, Y)$ is a global inverse semigroupoid action.

\begin{exe}
    Let $S = \{0,e,1\}$ be the inverse monoid from Example \ref{exe:0e1}, and let $(\theta, S)$ be the (in fact global) partial inverse semigroupoid action of $S$ on itself. Since $S$ is a semigroup and the action is non-degenerate, we have
        $$ D = S \times S \cup \{\delta\} \times S. $$
    The set $E$ obtained from our construction is $E = \{ [1,0], [1,e], [1,1], [0,1], [0,e], [e,1] \}$, where
    \begin{align*}
        [1,0] &= \{ (1,0), (e,0), (0,0), (\delta,0) \}, & [0,1] &= \{(0,1)\}, \\
        [1,e] &= \{ (1,e), (e,e), (\delta,e) \}, & [0,e] &= \{(0,e)\}, \\
        [1,1] &= \{ (1,1), (\delta,1) \}, & [e,1] &= \{ (0,1) \}.
    \end{align*}
    Since $S$ is a semigroup and the action is non-degenerate, we have ${_0E} = E = {_eE} = {_1E}$. The action $\beta_s$ is given by left multiplication by $s$ on the first coordinate of $[t,x]$. For instance, $\beta_0([1,1]) = [0,1] \neq [1,1]$. Therefore, $(\beta, E)$ fails to satisfy condition \eqref{I1}, and thus is not a partial inverse semigroupoid action.

    Now consider morphisms of partial actions from $(\beta, E)$ to $(\alpha, S)$.
    The only possible morphism $\Phi \colon (\beta, E) \to (\alpha, S)$ is $\Phi([t, x]) = 0$ for all $[t, x] \in E$. Indeed, since $E = {_0E}$, any morphism of partial actions must satisfy $\Phi(E) = \Phi({_0E}) \subseteq {_0S} = {0}$.
    
    On the other hand, $(\alpha, S)$ is a global inverse semigroupoid action, and the identity map $id_S \colon (\alpha, S) \to (\alpha, S)$ is itself a morphism of partial actions.
    However, $\Phi \circ \iota \equiv 0 \neq id_S$. Thus, $id_S$ cannot be factored through the embedding $\iota$.
\end{exe}

    
    \bibliographystyle{abbrvnat}

\end{document}